\documentclass[11pt]{article}
\usepackage{amssymb}
\makeatletter\@addtoreset {equation}{section}\makeatother

\newtheorem{theorem}{Theorem}
\newtheorem{lemma}{Lemma}
\newtheorem{remark}{Remark}
\newtheorem{corollary}{Corollary}
\newtheorem{assumption}{Assumption}
\newtheorem{example}{Example}
\newtheorem{definition}{Definition}

\newenvironment{proof}{
    \noindent {\it Proof.}}{\hfill$\Box$
}
\newenvironment{proof1}{
    \noindent {\it Proof }}{\hfill$\Box$
}

\usepackage[dvips]{epsfig}
\usepackage{graphicx}

\begin{document}

\title{\bf Moving gap solitons in periodic potentials}

\author{Dmitry Pelinovsky\footnote{\small On leave from Department of Mathematics, McMaster
University, Hamilton, Ontario, Canada, L8S 4K1} and Guido Schneider \\
{\small Institut f\"{u}r Analysis, Dynamik und
Modellierung Fakult\"{a}t f\"{u}r Mathematik und Physik,} \\
{\small Universit\"{a}t Stuttgart, Pfaffenwaldring 57, D-70569
Stuttgart, Germany} }

\date{\today}
\maketitle

\begin{abstract}
We address existence of moving gap solitons (traveling localized
solutions) in the Gross-Pitaevskii equation with a small periodic
potential. Moving gap solitons are approximated by the explicit
localized solutions of the coupled-mode system. We show however
that exponentially decaying traveling solutions of the
Gross-Pitaevskii equation do not generally exist in the presence
of a periodic potential due to bounded oscillatory tails ahead and
behind the moving solitary waves. The oscillatory tails are not
accounted in the coupled-mode formalism and are estimated by using
techniques of spatial dynamics and local center-stable manifold
reductions. Existence of bounded traveling solutions of the
Gross--Pitaevskii equation with a single bump surrounded by
oscillatory tails on a finite large interval of the spatial scale
is proven by using these technique.  We also show generality of
oscillatory tails in other nonlinear equations with a periodic
potential.
\end{abstract}

\section{Introduction}

Moving gap solitons are thought to be steadily traveling localized
solutions of nonlinear partial differential equations with
spatially periodic coefficients. The name of gap solitons comes
from the fact that parameters of {\em stationary} localized
solutions reside in the spectral gap of the associated
Schr\"{o}dinger operator. Existence of stationary solutions can be
studied by separation of the time and space variables and
reduction of the problem to an elliptic semi-linear equation.
Since the variables are not separable for {\em traveling}
solutions, a little is known about existence of moving gap
solitons.

We address stationary and traveling localized solutions in the
context of the Gross-Pitaevskii equation with an external periodic
potential
\begin{equation}
\label{GP} i E_t = - E_{xx} + V(x) E + \sigma |E|^2 E,
\end{equation}
where $E(x,t) : \mathbb{R} \times \mathbb{R}_+ \mapsto
\mathbb{C}$, $V(x) : \mathbb{R} \mapsto \mathbb{R}$, and $\sigma =
\pm 1$. The Gross--Pitaevskii equation (\ref{GP}) is derived for
the mean-field amplitude of the Bose--Einstein condensate placed
in the optical lattice trap $V(x)$, where $\sigma$ is the
normalized scattering length \cite{bec_book}. Stationary solutions
of the Gross--Pitaevskii equation are found from the solutions of
the nonlinear differential equation
\begin{equation}
\label{stationary} \phi''(x) + \omega \phi(x) = V(x) \phi(x) +
\sigma |\phi(x)|^2 \phi(x),
\end{equation}
where $\phi(x) : \mathbb{R} \mapsto \mathbb{C}$, $\omega \in
\mathbb{R}$, and the exact reduction $E(x,t) = \phi(x) e^{-i
\omega t}$ is used. Localized stationary solutions of the ODE
problem (\ref{stationary}) were proved to exist in \cite{Pankov}.

\begin{theorem}[{\bf Pankov, 05}] Let $V(x)$ be a real-valued,
bounded, and periodic potential. Let $\omega$ be in a finite gap
of the purely continuous spectrum of $L = -\partial_x^2 + V(x)$ in
$L^2(\mathbb{R})$. There exists a non-trivial weak solution
$\phi(x)$ in $H^1(\mathbb{R})$, which is (i) real-valued, (ii)
continuous on $x \in \mathbb{R}$ and (iii) decays exponentially as
$|x| \to \infty$.
\end{theorem}

In \cite{PSn07}, we have obtained a more precise information on
properties of the stationary solution $\phi(x)$ by working with a
small potential $V(x)$, when the spectrum of $L$ exhibits a
sequence of narrow gaps bifurcating near {\em resonance} points
$\omega = \omega_n = \frac{n^2}{4}$, $n \in \mathbb{N}$. We have
justified the use of the stationary coupled-mode equations which
have been used in the physics literature \cite{SS} for explicit
approximations of stationary gap solitons.

In this paper, we shall investigate whether the time-dependent
coupled-mode equations can be used for approximation of moving gap
solitons in the framework of the Gross--Pitaevskii equation
(\ref{GP}). The coupled-mode equations are derived formally when
the potential is represented by $V = \epsilon W(x)$, where
$\epsilon$ is small parameter and $W(x)$ is a normalized
$\epsilon$-independent function described as follows.

\begin{assumption}
Let $W(x)$ be a smooth $2\pi$-periodic function with zero mean and
symmetry $W(-x) = W(x)$ on $x \in \mathbb{R}$. The Fourier series
representation of $W(x)$ is
\begin{equation}
W(x) = \sum\limits_{m \in \mathbb{Z}} w_{2m} e^{i m x}, \quad
\mbox{such that} \quad \sum_{m \in \mathbb{Z}} (1 + m^2)^s
|w_{2m}|^2 < \infty, \;\; \forall s > \frac{1}{2},
\end{equation}
where $w_0 = 0$ and $w_{2m} = w_{-2m} = \bar{w}_{2m}$, $\forall m
\in \mathbb{N}$. \label{assumption-potential}
\end{assumption}

An asymptotic solution of the Gross--Pitaevskii equation
(\ref{GP}) is represented in the form \cite{GWH,SU}:
\begin{equation}
\label{formal_expansion} E(x,t) = \epsilon^{1/2} \left[ a(\epsilon
x,\epsilon t) e^{\frac{i n x}{2}} + b(\epsilon x,\epsilon t)
e^{-\frac{i n x}{2}} \right] e^{-\frac{i n^2 t}{4}} + {\rm
O}(\epsilon^{3/2}),
\end{equation}
where the vector function $(a,b) : \mathbb{R} \times \mathbb{R}_+
\mapsto \mathbb{C}^2$ satisfies the coupled-mode system
\begin{equation}
\label{cme} \left\{ \begin{array}{ccc} i ( a_T + n a_X) & = &
w_{2n} b + \sigma (|a|^2 + 2 |b|^2) a, \\ i ( b_T - n b_X) & = &
w_{2n} a + \sigma (2 |a|^2 + |b|^2) b,
\end{array} \right.
\end{equation}
in slow variables $X = \epsilon x$ and $T = \epsilon t$. The
system (\ref{cme}) admits a separation of variables \cite{SS}:
\begin{equation}
\label{traveling-solitons} a = \left( \frac{n+c}{n-c}\right)^{1/4}
A(\xi) e^{- i \mu n \tau}, \quad  b = \left( \frac{n-c}{n+c}
\right)^{1/4} B(\xi) e^{- i \mu n \tau},
\end{equation}
where $|c| < n$, the new independent coordinates $(\xi,\tau)$ are
given by the Lorentz transformation
\begin{equation}
\label{Lorentz} \xi = \frac{X - c T}{\sqrt{n^2 - c^2}}, \quad \tau
= \frac{T - c X}{\sqrt{n^2 - c^2}},
\end{equation}
and the new functions $A(\xi)$ and $B(\xi)$ satisfy the
coupled-mode system
\begin{equation}
\label{cme-ODE} \left\{ \begin{array}{ccc} (n-c) \left( i A' -
w_{2n} B \right) + \mu n (1-cn) A & = & \sigma \left[
(n+c) |A|^2 + 2 (n-c) |B|^2 \right] A, \\
-(n+c) \left( i B' + w_{2n} A \right) + \mu n (1 + cn) B & = &
\sigma \left[ 2 (n+c) |A|^2 + (n-c) |B|^2 \right] B.
\end{array} \right.
\end{equation}
Since $|A|^2 - |B|^2$ is constant in $\xi \in \mathbb{R}$ and the
constant is zero for localized solutions at infinity $|\xi| \to
\infty$, we can further represent the localized solution in the
form
\begin{equation}
\label{A-B-phi} A = \phi(\xi) e^{i \varphi(\xi)}, \qquad B =
\bar{\phi}(\xi) e^{i \varphi(\xi)},
\end{equation}
where the functions $\phi : \mathbb{R} \mapsto \mathbb{C}$ and
$\varphi : \mathbb{R} \mapsto \mathbb{R}$ are solutions of the
first-order equations
\begin{equation}
\label{varphi-phi} \left\{ \begin{array}{ccl} \varphi' & = &
\frac{n c \left( \mu (1 - n^2) - 2 \sigma |\phi|^2 \right)}{(n^2 -
c^2)}, \\ i  \phi' & = & w_{2n} \bar{\phi} - \frac{\mu n^2 (1 -
c^2)}{(n^2 - c^2)} \phi + \sigma \frac{(3 n^2 - c^2)}{(n^2 - c^2)}
|\phi|^2 \phi. \end{array} \right.
\end{equation}
The second equation of the system (\ref{varphi-phi}) is closed on
$\phi(\xi)$ and the explicit localized solution for $c \neq 0$ can
be found from the corresponding solution for $c = 0$ \cite{PSn07}.
For instance, if $\sigma = -1$ and $w_{2n} > 0$, the function
$\phi(\xi)$ is found in the explicit form
\begin{equation}
\phi = \sqrt{\frac{2(n^2-c^2)}{(3 n^2 - c^2)}}
\frac{\sqrt{w_{2n}^2 - \mu_0^2}}{\sqrt{w_{2n} + \mu_0} \; \cosh
\left( \sqrt{w_{2n}^2 - \mu_0^2} \xi \right) - i \sqrt{w_{2n} -
\mu_0} \; \sinh \left( \sqrt{w_{2n}^2 - \mu_0^2} \xi
\right)},\label{gap_soliton}
\end{equation}
where $\mu_0 = \mu \frac{n^2 (1 - c^2)}{(n^2 - c^2)}$ and $|\mu_0|
< w_{2n}$. In the case $c = 0$, we have $\mu_0 = \mu$ and the
condition $|\mu| < w_{2n}$ indicates that the frequency parameter
$\omega = \omega_n + \epsilon \mu$ of the stationary gap soliton
with $c = 0$ is chosen inside the newly formed gap of the
continuous spectrum near the bifurcation point $\omega_n =
\frac{n^2}{4}$, $n \in \mathbb{N}$ \cite{PSn07}, such that
\begin{equation}
\omega_n -\epsilon w_{2n} < \omega < \omega_n + \epsilon w_{2n}.
\end{equation}
The exact solution (\ref{gap_soliton}) can be extended easily to
values $\sigma = -1$ and $w_{2n} < 0$. Given a localized solution
for $\phi(\xi)$, we can integrate the first equation of the system
(\ref{varphi-phi}) and obtain a linearly growing solution for
$\varphi(\xi)$:
\begin{equation}
\label{phi-explicit} \varphi = \frac{n c}{(n^2 - c^2)} \left( \mu
(1 - n^2) \xi - 2 \sigma \int_0^{\xi} |\phi(\xi')|^2 d\xi'
\right).
\end{equation}
The trivial parameters of translations of solutions in $\xi$ and
$\varphi$ are set to zero in the explicit solutions
(\ref{gap_soliton}) and (\ref{phi-explicit}), such that the
functions $A(\xi)$ and $B(\xi)$ given by the parametrization
(\ref{A-B-phi}) satisfy the constraints $A(\xi) = \bar{A}(-\xi)$
and $B(\xi) = \bar{B}(-\xi)$.

\begin{definition}
\label{definition-orbit} The traveling solution of the
coupled-mode system (\ref{cme}) is said to be a {\em reversible
homoclinic orbit} if it decays to zero at infinity and satisfies
the constraints $A(\xi) = \bar{A}(-\xi)$ and $B(\xi) =
\bar{B}(-\xi)$ in parametrization
(\ref{traveling-solitons})--(\ref{A-B-phi}).
\end{definition}

We study persistence of the traveling solution
(\ref{traveling-solitons})--(\ref{phi-explicit}) with $c \neq 0$
of the coupled-mode system (\ref{cme}) in the Gross--Pitaevskii
equation (\ref{GP}). We show that the moving gap solitons have
bounded oscillatory tails in the far-field profile of the scale
$\epsilon^{N+1}$, which are small in amplitude of the order
$\epsilon^{N + 1/2}$ for any $N \geq 1$. These small oscillatory
tails are not accounted in the coupled-mode system (\ref{cme}).
The main result is formulated below.

\begin{theorem}
\label{theorem-main} Let Assumption \ref{assumption-potential} be
satisfied. Fix $n \in \mathbb{N}$, such that $w_{2n} \neq 0$. Let
$\omega = \frac{n^2}{4} + \epsilon \Omega$, such that $|\Omega| <
\Omega_0 = |w_{2n}| \frac{\sqrt{n^2-c^2}}{n}$. Let $0 < c < n$,
such that $\frac{n^2+c^2}{2c} \notin \mathbb{Z}'$, where
$\mathbb{Z}'$ is a set of odd (even) numbers for odd (even) $n$.
Fix $N \in \mathbb{N}$. For sufficiently small $\epsilon$, there
are $\epsilon$-independent constants $L > 0$ and $C > 0$, such
that there exists an infinite-dimensional, continuous family of
traveling solutions of the Gross--Pitaevskii equation (\ref{GP})
in the form $E(x,t) = e^{-i \omega t} \psi(x,y)$, where $y = x -
ct$ and the function $\psi(x,y)$ is periodic (anti-periodic)
function of $x$ for even (odd) $n$, satisfying the reversibility
constraint $\psi(x,y) = \bar{\psi}(x,-y)$, and the bound
\begin{equation}
\label{bound-main} \left| \psi(x,y) - \epsilon^{1/2} \left[
a_{\epsilon} (\epsilon y) e^{\frac{i n x}{2}} +
b_{\epsilon}(\epsilon y) e^{-\frac{i n x}{2}}\right] \right| \leq
C_0 \epsilon^{N+1/2}, \quad \forall x \in \mathbb{R}, \; \forall y
\in [-L/\epsilon^{N+1},L/\epsilon^{N+1}],
\end{equation}
Here $a_{\epsilon}(Y) = a(Y) + {\rm O}(\epsilon)$ and
$b_{\epsilon}(Y) = b(Y) + {\rm O}(\epsilon)$ on $Y = \epsilon y
\in \mathbb{R}$ are exponentially decaying solutions as $|Y| \to
\infty$, where $a(Y)$ and $b(Y)$ are solutions of the coupled-mode
system (\ref{cme}) with $Y = X - cT$.
\end{theorem}

\begin{remark}
{\rm      $\phantom{t}$

\noindent (a) The solution $\psi(x,y)$ is a bounded non-decaying
function on a large finite interval $$y \in [-L/\epsilon^{N+1},
L/\epsilon^{N+1}] \subset \mathbb{R}$$ but we do not claim that
the solution $\psi(x,y)$ can be extended to a global bounded
function on $y \in \mathbb{R}$.

\noindent (b) Since the homoclinic orbit $(a,b)$ of the
coupled-mode system (\ref{cme}) is single-humped, the traveling
solution $\psi(x,y)$ is represented by a single bump surrounded by
bounded oscillatory tails.

\noindent (c) The solution $(a_{\epsilon},b_{\epsilon})$ is
defined up to the terms of ${\rm O}(\epsilon^N)$ and it satisfies
an extended coupled-mode system with a unique reversible
single-humped homoclinic orbit. }
\end{remark}

Our work can be compared with three groups of papers. The first
group covers rigorous justification of the validity of the
coupled-mode system (\ref{cme}) for the system of cubic Maxwell
equations \cite{GWH} and for the Klein--Fock equation with
quadratic nonlinearity \cite{SU}. The bound on the error terms was
proved for a finite time interval, which depends on $\epsilon$. By
using this bound, one can see that the solution of the
Gross--Pitaevskii equation (\ref{GP}) behaves as a moving gap
soliton of the coupled-mode system (\ref{cme}) during the initial
time evolution in $H^1(\mathbb{R})$ \cite{GWH} or in
$C_b^0(\mathbb{R})$ \cite{SU}. However, the error is not
controlled on the entire time interval $t \in \mathbb{R}$ since
other effects such as radiation due to interactions of the moving
gap soliton with the stationary periodic potential can destroy
steady propagation of gap solitons.

The second group of articles covers analysis of persistence of
small-amplitude localized modulated pulses in nonlinear dispersive
systems such as the Maxwell equations with periodic coefficients
\cite{AK}, the nonlinear wave equation \cite{GS01}, and the
quasilinear wave equation \cite{GS05}. Methods of spatial dynamics
were applied in these works to show that a local center manifold
spanned by oscillatory modes destroys exponential localization of
the modulated pulses along the directions of the slow stable and
unstable manifolds. As a result, the modulating pulse solutions
decay in the spatial dynamics to small-amplitude oscillatory
disturbances in the far-field regions.

The third group of papers addresses propagation of a moving
solitary wave in a periodic potential $V(x)$ of a large period
(see review in \cite{SB}). An effective particle equation is
derived from the focusing Gross--Pitaevskii equation (\ref{GP})
with $\sigma = -1$ by an heuristic asymptotic expansion. The
particle equation describes a steady propagation of the moving
solitary wave with $\omega < 0$, which corresponds to the
semi-infinite gap of the periodic potential. Radiation effects
appear beyond all orders of the asymptotic expansion. They have
been incorporated in the asymptotic formalism by using
perturbation theory based on the inverse scattering transform
\cite{KM}. The same methods were applied to the finite-period and
small-period potentials \cite{SB}. Unfortunately, this group of
article does not connect individual results in a complete rigorous
theory of the time evolution of a solitary wave in a periodic
potential, although it does gives a good intuition on what to
expect from the time evolution.

Our article is structured as follows. Section 2 reformulates the
existence problem for moving gap solitons as the spatial dynamical
system. Section 3 presents the Hamiltonian structure for the
spatial dynamical system and normal coordinates of the Hamiltonian
system. Section 4 describes a transformation of the Hamiltonian
system to the normal form and gives a proof of persistence of a
reversible homoclinic orbit in the extended coupled-mode system.
Section 5 presents a construction of a local center-saddle
manifold which concludes the proof of Theorem \ref{theorem-main}.
Section 6 discusses other models for moving gap solitons with
oscillatory tails.

\section{Spatial dynamics formulation}

We look for traveling solutions of the Gross--Pitaevskii equation
(\ref{GP}) in the form
\begin{equation}
E(x,t) = e^{-i \omega t} \psi(x,y), \qquad y = x - ct,
\end{equation}
where $\omega$ is a parameter of gap solitons and the coordinates
$(x,y)$ are linearly independent if $c \neq 0$. For simplicity, we
only consider the case $c > 0$. The envelope function $\psi(x,y)$
satisfies the partial differential equation
\begin{equation}
\label{PDE} \left( \omega - i c \partial_y + \partial_x^2 + 2
\partial_x \partial_y + \partial_y^2 \right) \psi(x,y) =
\epsilon W(x) \psi(x,y) + \sigma |\psi(x,y)|^2 \psi(x,y).
\end{equation}
At this stage, the equation (\ref{PDE}) is {\em equivalent} to the
original equation (\ref{GP}) if $c \neq 0$. We shall however
specify the class of functions $\psi(x,y)$ to accommodate the
moving gap solitons according to their leading-order
representation given by (\ref{formal_expansion}),
(\ref{traveling-solitons}), and (\ref{Lorentz}). In particular, we
consider either periodic (for even $n$) or anti-periodic (for odd
$n$) functions $\psi(x,y)$ in variable $x \in [0,2\pi]$ and look
for decaying or bounded solutions $\psi(x,y)$ in variable $y \in
\mathbb{R}$. Such solutions can be described by using the
formalism of spatial dynamical systems \cite{GS01,GS05}. We make
use of the periodic or anti-periodic conditions in variable $x$
and represent the solution $\psi(x,y)$ in the form
\begin{equation}
\label{solution-series} \psi(x,y) = \sqrt{\epsilon} \sum_{m \in
\mathbb{Z}'} \psi_m(y) e^{\frac{i}{2} m x}, \quad \psi_m(y) =
\frac{1}{2 \pi \sqrt{\epsilon}} \int_0^{2 \pi} \psi(x,y)
e^{-\frac{i}{2} m x} dx, \quad m \in \mathbb{Z}',
\end{equation}
where the factor $\sqrt{\epsilon}$ is used for the convenience and
the set $\mathbb{Z}'$ contains even numbers if $\psi(x,y)$ is
periodic in $x$ and odd numbers if $\psi(x,y)$ is anti-periodic in
$x$. The series representation (\ref{solution-series}) transforms
the PDE system (\ref{PDE}) to the nonlinear system of coupled ODEs
\begin{eqnarray}
\nonumber \psi_m''(y) + i (m - c) \psi_m'(y) + \left(\omega -
\frac{m^2}{4} \right) \psi_m(y) = \epsilon \sum_{m_1 \in
\mathbb{Z}'} w_{m-m_1} \psi_{m_1}(y) \\
\label{difference-system}
 + \epsilon \sigma \sum_{m_1 \in
\mathbb{Z}'}\sum_{m_2 \in \mathbb{Z}'} \psi_{m_1}(y)
\bar{\psi}_{-m_2}(y) \psi_{m - m_1 - m_2}(y), \quad m \in
\mathbb{Z}'.
\end{eqnarray}
The left-hand-side of the system (\ref{difference-system})
represents a linearized system at the zero solution for $\epsilon
= 0$. Since the linearized system at $\epsilon = 0$ has a diagonal
structure on $m \in \mathbb{Z}'$, its solutions are given by the
eigenmodes $\psi_{m'}(y) = e^{\kappa_m y} \delta_{m,m'}$ with
$m,m' \in \mathbb{Z}'$, where the values of $\kappa_m$ are
determined by the roots of quadratic equations
\begin{equation}
\label{quadratic-equations} \kappa = \kappa_m : \quad \kappa^2 +
i(m-c) \kappa + \omega - \frac{m^2}{4} = 0, \qquad m \in
\mathbb{Z}'.
\end{equation}
The zero root $\kappa = 0$ exists if and only if $\omega =
\omega_n = \frac{n^2}{4}$ for any fixed $n \in \mathbb{Z}$. The
zero root has multiplicity two for $m = \pm n$ if $n \in
\mathbb{N}$ and $c \neq \pm n$. The special value $\omega =
\omega_n$ corresponds to the bifurcation of periodic or
anti-periodic solutions as well as of the stationary gap solitons
with $c = 0$ \cite{PSn07}. We note that the values of $n$
determine the choice for the set $\mathbb{Z}'$: it includes even
(odd) numbers if $n$ is even (odd). We shall hence focus on the
bifurcation case $\omega = \omega_n$, when the two roots of the
quadratic equations (\ref{quadratic-equations}) are represented
explicitly as follows
\begin{equation}
\label{root-quadratic} \omega = \omega_n : \quad \kappa =
\kappa^{\pm}_m = \frac{i(c-m) \pm \sqrt{2c m - n^2-c^2}}{2},
\qquad m \in \mathbb{Z}'.
\end{equation}
When $m > m_0 = \left[ \frac{n^2+c^2}{2c} \right]$, where the
notation $[a]'$ denotes the integer part of the number $a$ in the
set $\mathbb{Z}'$, all roots are complex-valued with ${\rm
Re}(\kappa^{\pm}_m) = \pm \frac{1}{2} \sqrt{2c m - n^2 - c^2} \neq
0$ and ${\rm Im}(\kappa^{\pm}_m) = \frac{c-m}{2}$. When $m \leq
m_0$, al roots $\kappa$ are purely imaginary with $\kappa_m^{\pm}
= i k_m^{\pm}$ and $k_m^{\pm} = \frac{(c-m) \pm \sqrt{n^2+c^2-2c
m}}{2}$.

\begin{lemma}
\label{lemma-roots} Let $\omega = \omega_n$, $n \in \mathbb{N}$
and $c > 0$, such that $\frac{n^2+c^2}{2c} \notin \mathbb{Z}'$.
Then,

\begin{itemize}
\item[(i)] The phase space of the linearized system
(\ref{difference-system}) at the zero solution for $\epsilon = 0$
decomposes into a direct sum of subspaces $E^s \oplus E^u \oplus
E^{c^+} \oplus E^{c^-}$, where
\begin{equation}
\label{subspaces} E^s = \oplus_{m > m_0} E_m^+, \quad E^u =
\oplus_{m > m_0} E_m^-, \quad E^{c^+} = \oplus_{m \leq m_0} E_m^+,
\quad E^{c^-} = \oplus_{m \leq m_0} E_m^-.
\end{equation}
The subspace $E_m^{\pm}$ consists of the eigenspace associated
with the $m$-th Fourier component of the solution
(\ref{solution-series}) corresponding to the root $\kappa =
\kappa^{\pm}_m$ in (\ref{root-quadratic}).

\item[(ii)] The zero root $\kappa = 0$ is semi-simple of
multiplicity two. The purely imaginary roots $\kappa \in i
\mathbb{R}$ are semi-simple of the maximal multiplicity three. All
other roots $\kappa \in \mathbb{C}$ are simple.
\end{itemize}
\end{lemma}

\begin{proof}
It follows from the quadratic equation (\ref{quadratic-equations})
that a root $\kappa = \kappa_0$ is double if $\kappa_0 =
\frac{i(c-m)}{2}$, which implies that $m = \frac{n^2+c^2}{2c}$.
Under the non-degeneracy constraint $\frac{n^2+c^2}{2c} \notin
\mathbb{Z}'$, all roots $\kappa$ are semi-simple. When $m > m_0 =
\left[\frac{n^2+c^2}{2c} \right]'$, all roots are complex-valued
and simple, such that $E^s$ and $E^u$ are stable and unstable
manifolds of the linearized system at the zero solution for
$\epsilon = 0$.

When $m \leq m_0$, all roots $\kappa$ are purely imaginary, such
that $E^{c^+} \oplus E^{c^-}$ is a center manifold of the
linearized system for $\epsilon = 0$. It is obvious that
$\kappa^+_m$ increases as $m$ decreases, while $\kappa^-_m$
decreases for $m_1 < m \leq m_0$ and increases for $m \leq m_1$ as
$m$ decreases, where $m_1 = \left[ \frac{n^2}{2c} \right]'$.
Therefore, the purely imaginary roots may have the maximal
multiplicity three. The zero eigenvalue has however multiplicity
two since the two modes $m = n$ and $m = - n$ have simple zero
eigenvalues and other modes have no zero eigenvalues.
\end{proof}

\begin{lemma}
Let $\omega = \omega_n$, $n \in \mathbb{N}$ and $0 < c < n$. If
$c$ is irrational, all non-zero roots $\kappa$ of the quadratic
equations (\ref{quadratic-equations}) are simple.
\label{lemma-multiplicity}
\end{lemma}

\begin{proof}
By Lemma \ref{lemma-roots}, only imaginary roots $\kappa$ can be
semi-simple. Let two roots $\kappa$ coincide for $m \leq m_0$ and
$l \leq m_0$. Then, $\kappa = -i \frac{m+l}{4}$ and $(m,l)$
satisfies the equation
\begin{equation}
\label{semisimple-equation} (m-l)^2 + 4 c (m+l) - 4 n^2 = 0,
\qquad m \leq m_0, \; l \leq m_0.
\end{equation}
If $c$ is irrational, equation (\ref{semisimple-equation}) has no
solutions for integers $m$ and $l$. Therefore, all non-zero roots
$\kappa$ are simple.
\end{proof}

\begin{example}
{\rm Figure 1(a) illustrates the distribution of imaginary roots
$\kappa^{\pm}_m = i k^{\pm}_m$, $m \leq m_0$ for $c =
\frac{1}{\sqrt{2}}$ and $n = 1$. Although all imaginary roots are
simple for this (irrational) value of $c$, the purely imaginary
roots $\kappa$ can approach to each other arbitrarily close.
Figure 1(b) shows a similar distribution for $c = \frac{1}{2}$ and
$n = 1$. It follows from equation (\ref{semisimple-equation}) that
an infinite sequence of semi-simple roots exists for this
(rational) value of $c$.}
\end{example}

\begin{figure}[htbp]
\begin{center}
\includegraphics[height=6cm]{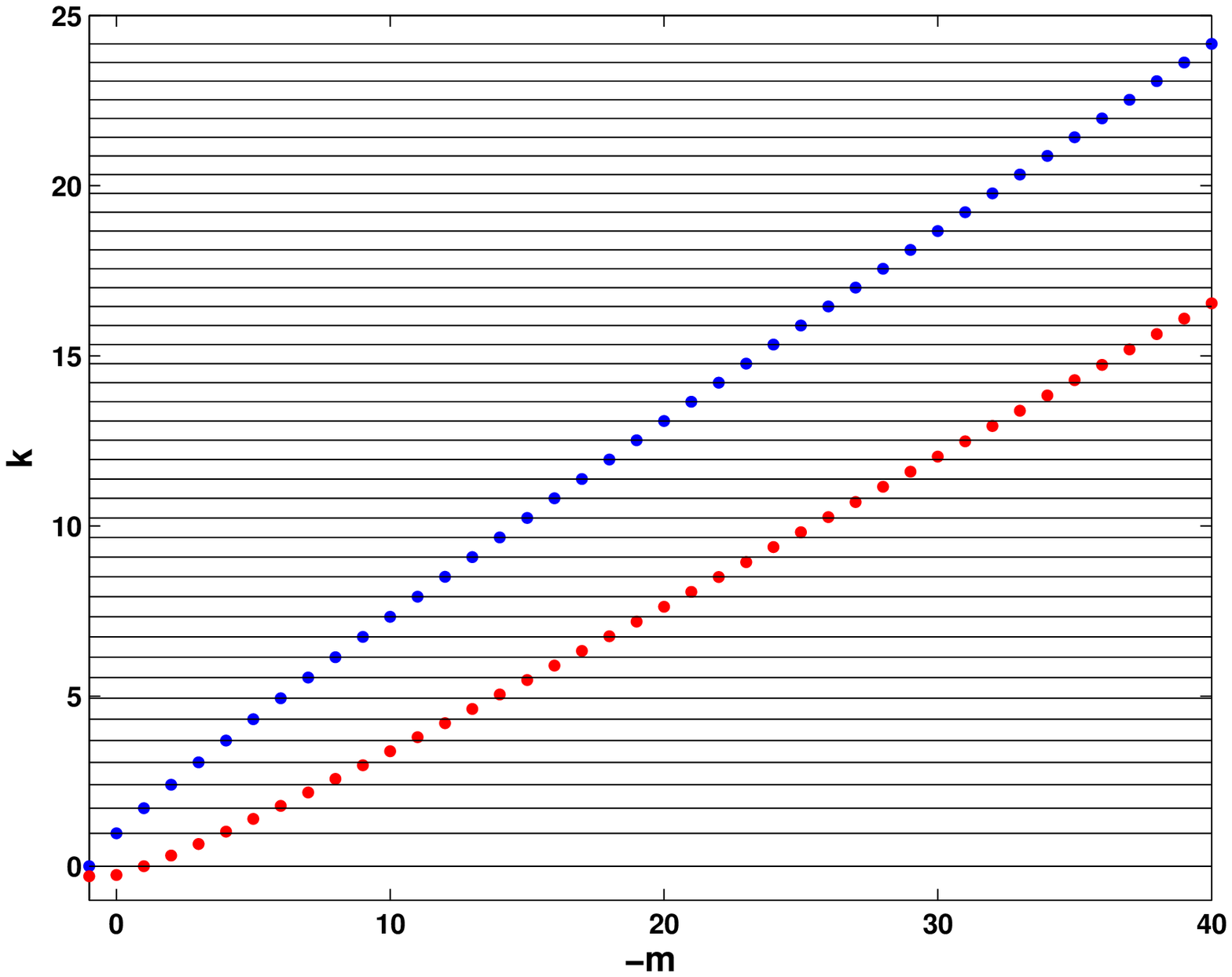}
\includegraphics[height=6cm]{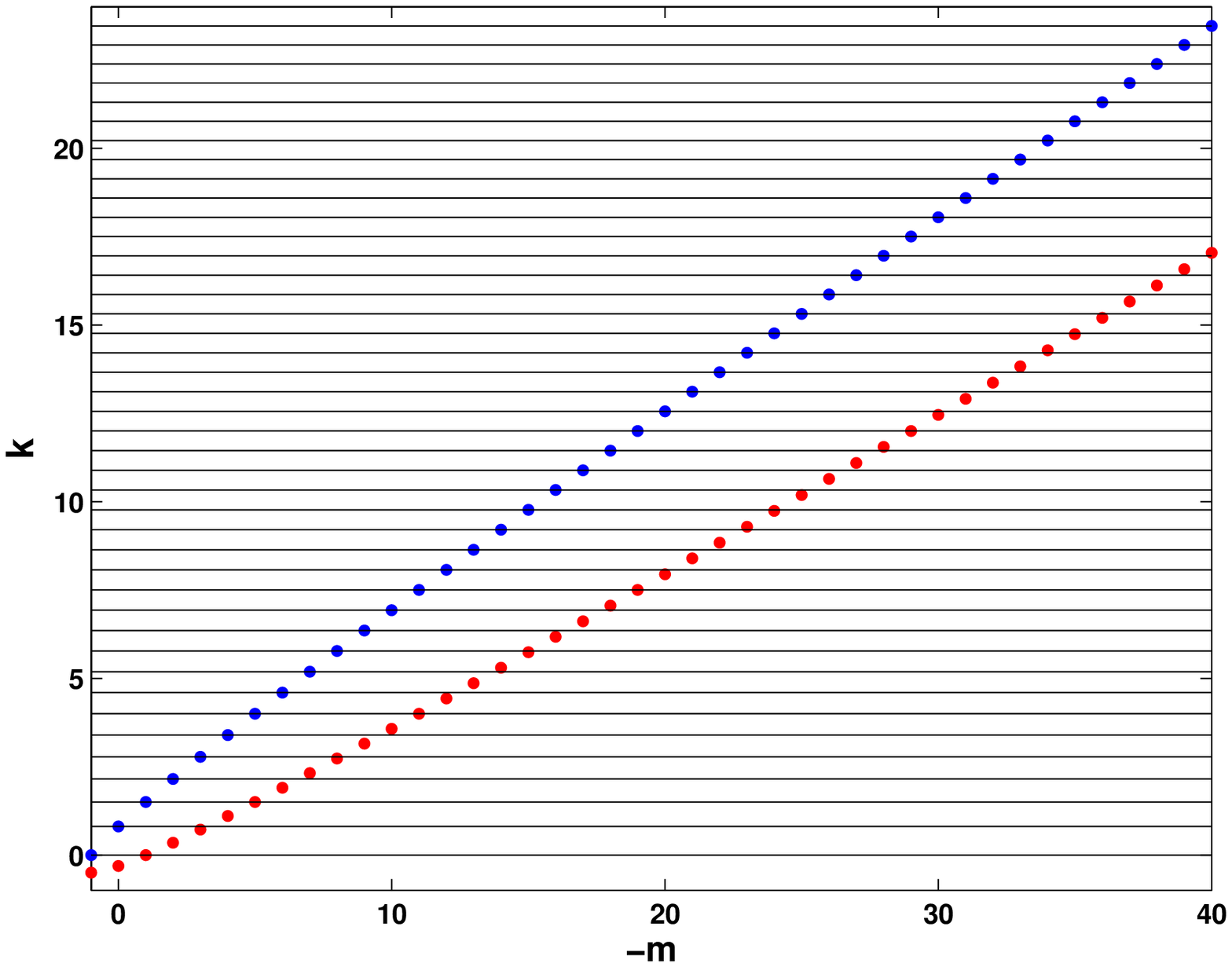}
\end{center}
\caption{Purely imaginary roots $\kappa^{\pm}_m = i k_m^{\pm}$
versus $m \leq m_0$ for $n = 1$, $c = \frac{1}{\sqrt{2}}$ (left)
and for $n = 1$, $c = \frac{1}{2}$ (right). Blue dots correspond
to $k_m^+$ and red dots correspond to $k_m^-$.} \label{fig-roots}
\end{figure}

\begin{remark}
{\rm When $c = 0$, the formalism of spatial dynamical systems
fails. Indeed, the linearized system (\ref{difference-system}) for
$\epsilon = 0$ and $c = 0$ has a set of semi-simple eigenvalues of
multiplicity two at $\kappa = \pm i k$, $k \in \mathbb{Z}$.
However, all non-zero eigenvalues are equivalent to the zero
eigenvalue due to existence of infinitely many symmetries for the
system (\ref{difference-system}): $\psi_m(y) \mapsto
\psi_{m+m'}(y) e^{\frac{i m' y}{2}}$, $\forall m' \in
\mathbb{Z}'$. Therefore, the behavior of all modes with $m \in
\mathbb{Z}' \backslash \{ n,-n\}$ repeat the behavior of the two
resonant modes with $m = \pm n$ and only these two modes are
relevant for existence of localized solutions for $\epsilon \neq
0$ and $c = 0$. The failure of the spatial dynamics formulation is
related to the fact that $y = x$ for $c = 0$, such that the
second-order ODE (\ref{stationary}) can not be replaced by the PDE
(\ref{PDE}) and hence it can not be written as a system of
infinitely many second-order ODEs (\ref{difference-system}). }
\end{remark}

\begin{lemma}
\label{lemma-inversion} Consider a linear inhomogeneous equation
\begin{equation}
\label{linear-inhomogeneous} \left( \partial_y^2 + i (m-c)
\partial_y + \frac{n^2-m^2}{4} \right) \psi_m(y) = -F_m(y), \qquad
\forall m > m_0,
\end{equation}
where $F_m \in C_b^0(\mathbb{R})$. There exists a unique solution
$\psi_m \in C_b^0(\mathbb{R})$, such that $\| \psi_m
\|_{C_b^0(\mathbb{R})} \leq C \| F_m \|_{C_b^0(\mathbb{R})}$ for
some $C > 0$.
\end{lemma}

\begin{proof}
Let $\psi_m = e^{-i \alpha y} \varphi_m(y)$ with $\alpha =
\frac{m-c}{2}$ and rewrite the linear equation
(\ref{linear-inhomogeneous}) in the equivalent form
$$
\left( \beta^2 - \partial_y^2 \right) \varphi_m(y) = F_m(y) e^{i
\alpha y}, \qquad \beta^2 = \frac{2 c m - n^2 - c^2}{4} > 0.
$$
Since solutions of the homogeneous equation are exponentially
decaying and growing as $\varphi_m \sim e^{\pm \beta y}$, there
exists a unique bounded solution of the inhomogeneous equation in
the integral form
$$
\varphi_m(y) = \frac{1}{2\beta} \int_{-\infty}^{\infty} e^{-\beta
|y - y'|} F_m(y') e^{i \alpha y'} dy',
$$
such that $\| \varphi_m \|_{C_b^0(\mathbb{R})} \leq
\frac{1}{\beta^2} \| F_m \|_{C_b^0(\mathbb{R})}$.
\end{proof}

\begin{remark}
{\rm By using Lemma \ref{lemma-inversion} and the Implicit
Function Theorem in suitable vector spaces, one can solve
equations of the system (\ref{difference-system}) for $m
> m_0$ and parameterize all components $\psi_m(y)$ with $m
> m_0$ by bounded components $\psi_m(y)$ with $m \leq m_0$ for sufficiently small
$\epsilon$. However, we do not perform this elimination at this
stage, since we are going to rewrite the system
(\ref{difference-system}) as a Hamiltonian dynamical system and
use a formalism of near-identity transformations and normal forms,
which works easier if the symplectic structure of the Hamiltonian
system is local. }
\end{remark}

\section{Hamiltonian formalism and normal coordinates}

We rewrite the system of second-order equations
(\ref{difference-system}) as the system of first-order equations
which admits a symplectic Hamiltonian structure. Let $\omega =
\frac{n^2}{4} + \epsilon \Omega$, where $n \in \mathbb{N}$ and
$\Omega$ is a free parameter. Let $\phi_m(y) = \psi_m'(y) -
\frac{i}{2}(c-m) \psi_m(y)$ for all $m \in \mathbb{Z}'$. The
system (\ref{difference-system}) is equivalent to the first-order
system
\begin{eqnarray}
\label{difference-system-first-order} \left\{ \begin{array}{ccl}
\frac{d \psi_m}{d y} & = & \phi_m + \frac{i}{2} (c - m) \psi_m
\\ \frac{d \phi_m}{d y} & = & -\frac{1}{4} \left( n^2 + c^2 - 2 c m
\right) \psi_m + \frac{i}{2} (c-m) \phi_m - \epsilon \Omega \psi_m
+ \epsilon \sum_{m_1 \in \mathbb{Z}'} w_{m-m_1} \psi_{m_1}
\\ & \phantom{t} & \phantom{texttext} + \epsilon \sigma \sum_{m_1
\in \mathbb{Z}'}\sum_{m_2 \in \mathbb{Z}'} \psi_{m_1}
\bar{\psi}_{-m_2} \psi_{m - m_1 - m_2}.
\end{array} \right.
\end{eqnarray}
Let bolded symbol $\mbox{\boldmath $\psi$}$ denote a vector
consisting of elements of the set $\{ \psi_m \}_{m \in
\mathbb{Z}'}$. The variables $\{\mbox{\boldmath
$\psi$},\mbox{\boldmath $\phi$},\bar{\mbox{\boldmath
$\psi$}},\bar{\mbox{\boldmath $\phi$}}\}$ are canonical and the
system (\ref{difference-system-first-order}) is equivalent to the
Hamilton's equations of motion
\begin{equation}
\label{Ham-system} \frac{d \psi_m}{d y} = \frac{\partial
H}{\partial \bar{\phi}_m}, \quad  \frac{d \phi_m}{d y} =
-\frac{\partial H}{\partial \bar{\psi}_m}, \quad m \in
\mathbb{Z}',
\end{equation}
where $H = H(\mbox{\boldmath $\psi$},\mbox{\boldmath
$\phi$},\bar{\mbox{\boldmath $\psi$}},\bar{\mbox{\boldmath
$\phi$}})$ is the Hamiltonian function given by
\begin{eqnarray}
\nonumber H & = & \sum_{m \in \mathbb{Z}'} \left[ |\phi_m|^2 +
\frac{1}{4} (n^2 + c^2 - 2 c m) |\psi_m|^2 +
\frac{i}{2} (c-m) (\psi_m \bar{\phi}_m - \bar{\psi}_m \phi_m )
+ \epsilon \Omega |\psi_m|^2 \right] \\
\label{Hamiltonian} & \phantom{t} & - \epsilon \sum_{m \in
\mathbb{Z}'} \sum_{m_1 \in \mathbb{Z}'} w_{m-m_1} \psi_{m_1}
\bar{\psi}_m - \frac{\epsilon \sigma}{2} \sum_{m \in \mathbb{Z}'}
\sum_{m_1 \in \mathbb{Z}'}\sum_{m_2 \in \mathbb{Z}} \psi_{m_1}
\bar{\psi}_{-m_2} \psi_{m-m_1-m_2} \bar{\psi}_m.
\end{eqnarray}
Let us define the discrete weighted $l^2$-space by its norm
\begin{equation}
\label{discrete-Sobolev-norm} \forall {\bf u} \in
l^2_s(\mathbb{Z}') : \quad \| {\bf u} \|^2_{l^2_s} = \sum_{m \in
\mathbb{Z}'} (1 + m^2)^s |u_m|^2 < \infty.
\end{equation}
Since $l^2_s(\mathbb{Z}')$ is Banach algebra for $s >
\frac{1}{2}$, the convolution sums in the nonlinear system
(\ref{difference-system-first-order}) are closed if ${\bf W} \in
l^2_s(\mathbb{Z})$ (Assumption \ref{assumption-potential}) and
$\mbox{\boldmath $\psi$} \in l^2_s(\mathbb{Z}')$ for $s >
\frac{1}{2}$. Due to the unbounded linear part, the vector field
of the system (\ref{difference-system-first-order}) map a domain
in $D$ and to a range in $X$, where $D$ and $X$ are given by
\begin{eqnarray}
D = \left\{ (\mbox{\boldmath $\psi$},\mbox{\boldmath
$\phi$},\bar{\mbox{\boldmath $\psi$}},\bar{\mbox{\boldmath
$\phi$}}) \in l^2_{s+1}(\mathbb{Z}',\mathbb{C}^4) \right\}, \qquad
X = \left\{ (\mbox{\boldmath $\psi$},\mbox{\boldmath
$\phi$},\bar{\mbox{\boldmath $\psi$}},\bar{\mbox{\boldmath
$\phi$}}) \in l^2_s(\mathbb{Z}',\mathbb{C}^4) \right\},
\end{eqnarray}
for any $s > \frac{1}{2}$. We note that $D \subset X$ and that $X$
can be chosen as the phase space of the Hamiltonian system
(\ref{Ham-system}).

If $W(-x) = W(x)$, $\forall x \in \mathbb{R}$ (Assumption
\ref{assumption-potential}), then $w_{2m} = w_{-2m}$, $\forall m
\in \mathbb{Z}$. In this case, the Hamiltonian system
(\ref{difference-system-first-order}) is reversible and its
solutions are invariant under the transformation
\begin{equation}
\label{transformation-reversibility} \mbox{\boldmath $\psi$}(y)
\mapsto \bar{\mbox{\boldmath $\psi$}}(-y), \quad \mbox{\boldmath
$\phi$}(y) \mapsto -\bar{\mbox{\boldmath $\phi$}}(-y).
\end{equation}
In addition, the Hamiltonian function (\ref{Hamiltonian}) is
invariant with respect to the gauge transformation
\begin{equation}
\label{transformation-gauge} \mbox{\boldmath $\psi$}(y) \mapsto
e^{i \alpha} \mbox{\boldmath $\psi$}(y), \;\; \mbox{\boldmath
$\phi$}(y) \mapsto e^{i \alpha} \mbox{\boldmath $\phi$}(y), \qquad
\forall \alpha \in \mathbb{R}.
\end{equation}
Let us define the domain for reversible solutions by
\begin{equation}
\label{domain-reversibility} D_r = \left\{ (\mbox{\boldmath
$\psi$},\mbox{\boldmath $\phi$},\bar{\mbox{\boldmath
$\psi$}},\bar{\mbox{\boldmath $\phi$}}) \in D : \quad
\mbox{\boldmath $\psi$}(-y) = \bar{\mbox{\boldmath $\psi$}}(y),
\;\;  \mbox{\boldmath $\phi$}(-y) = -\bar{\mbox{\boldmath
$\phi$}}(y) \right\}.
\end{equation}
If a local solution $(\mbox{\boldmath $\psi$}(y),\mbox{\boldmath
$\phi$}(y))$ of the system (\ref{difference-system-first-order})
is constructed on $y \in \mathbb{R}_+$ and it intersects at $y =
0$ with the reversibility constraint
\begin{equation}
\label{constraint-reversibility} \Sigma_r = \left\{
(\mbox{\boldmath $\psi$},\mbox{\boldmath
$\phi$},\bar{\mbox{\boldmath $\psi$}},\bar{\mbox{\boldmath
$\phi$}}) \in D : \quad {\rm Im} \mbox{\boldmath $\psi$} = 0, \;\;
{\rm Re} \mbox{\boldmath $\phi$} = 0 \right\},
\end{equation}
then the solution $(\mbox{\boldmath $\psi$}(y),\mbox{\boldmath
$\phi$}(y))$ is extended to a global reversible solution in $D_r$
on $y \in \mathbb{R}$ by using the reversibility transformation
(\ref{transformation-reversibility}). The global reversible
solution does not have an arbitrary parameter $\alpha$ induced by
the gauge transformation (\ref{transformation-gauge}).

In order to construct non-trivial bounded solutions of the
Hamiltonian system (\ref{difference-system-first-order}), we shall
introduce normal coordinates for the infinite-dimensional stable,
unstable and center manifolds of the linearized Hamiltonian system
(\ref{difference-system-first-order}) at the zero solution for
$\epsilon = 0$ (Lemma \ref{lemma-roots}). Let $\mathbb{Z}'_- = \{
m \in \mathbb{Z}' : \; m \leq m_0 \}$ and $\mathbb{Z}'_+ = \{ m
\in \mathbb{Z}' : \; m > m_0 \}$. For the center manifold of the
linearized Hamiltonian system, we set
\begin{equation}
\label{normal-coordinates-1} m \in \mathbb{Z}'_- : \quad \psi_m =
\frac{c_m^+(y) + c_m^-(y)}{\sqrt[4]{n^2 + c^2 - 2 c m}}, \qquad
\phi_m = \frac{i}{2} \sqrt[4]{n^2 + c^2 - 2 cm} \left[ c_m^+(y) -
c_m^-(y) \right].
\end{equation}
For the stable and unstable manifolds of the linearized
Hamiltonian system, we set
\begin{equation}
\label{normal-coordinates-2} m \in \mathbb{Z}'_+ : \quad \psi_m =
\frac{c_m^+(y) + c_m^-(y)}{\sqrt[4]{2 c m - n^2 - c^2}}, \qquad
\phi_m = \frac{1}{2} \sqrt[4]{2 c m - n^2 - c^2} \left[ c_m^+(y) -
c_m^-(y) \right].
\end{equation}
By using the normal coordinates
(\ref{normal-coordinates-1})--(\ref{normal-coordinates-2}), the
Hamiltonian function (\ref{Hamiltonian}) transforms to the new
form
\begin{eqnarray}
\label{Hamiltonian-transformed} H & = & \sum_{m \in \mathbb{Z}'_-}
\left( k_m^+ |c_m^+|^2 - k_m^- |c_m^-|^2 \right) + \sum_{m \in
\mathbb{Z}'_+} \left( \kappa_m^- c_m^- \bar{c}_m^+ - \kappa_m^+ c_m^+ \bar{c}_m^- \right) \\
\nonumber & \phantom{t} & + \epsilon \Omega \sum_{m \in
\mathbb{Z}'} \frac{|c_m^+ + c_m^-|^2}{\sqrt{|n^2 + c^2 - 2 cm|}} -
\epsilon \sum_{m \in \mathbb{Z}'} \sum_{m_1 \in \mathbb{Z}'}
w_{m,m_1} (c_{m_1}^+ + c_{m_1}^-) (\bar{c}_m^+ + \bar{c}_m^-) \\
\nonumber  & \phantom{t} & - \frac{\epsilon \sigma}{2} \sum_{m \in
\mathbb{Z}'} \sum_{m_1 \in \mathbb{Z}'}\sum_{m_2 \in \mathbb{Z}}
g_{m,m_1,m_2} (c_{m_1}^+ + c_{m_1}^-) (\bar{c}_{-m_2}^+ +
\bar{c}_{-m_2}^-)(c_{m-m_1-m_2}^+ + c_{m - m_1 - m_2}^-)
(\bar{c}_m^+ + \bar{c}_m^-),
\end{eqnarray}
where
\begin{eqnarray*}
& \phantom{t} & k_m^{\pm} = \frac{ c - m \pm \sqrt{n^2 + c^2 - 2
cm}}{2}, \qquad m \in \mathbb{Z}'_-, \\
& \phantom{t} & \kappa_m^{\pm} = \frac{i(c -
m) \pm \sqrt{2 cm - n^2 - c^2}}{2}, \qquad m \in \mathbb{Z}'_+, \\
& \phantom{t} & w_{m,m_1} = \frac{w_{m-m_1}}{\sqrt[4]{|n^2 + c^2 -
2 cm | | n^2 + c^2 - 2cm_1|}}, \qquad m,m_1 \in \mathbb{Z}'
\end{eqnarray*}
and
\begin{eqnarray*}
g_{m,m_1,m_2} = \frac{1}{\sqrt[4]{|n^2 + c^2 - 2 cm_1| |n^2 + c^2
+ 2 c m_2| |n^2 + c^2 - 2c(m-m_1-m_2)| |n^2 + c^2 - 2 c m|}},
\end{eqnarray*}
for all $m,m_1,m_2 \in \mathbb{Z}'$. The quadratic part of the
Hamiltonian function in (\ref{Hamiltonian-transformed}) for
$\epsilon = 0$ is diagonal in normal coordinates for $m \in
\mathbb{Z}_-'$ and it is block-diagonal for $m \in \mathbb{Z}_+'$.
The Hamiltonian equations of motions (\ref{Ham-system}) transform
in new canonical variables to the form
\begin{equation}
\label{Ham-system-transformed-1} \frac{d c_m^+}{d y} = i
\frac{\partial H}{\partial \bar{c}_m^+}, \quad  \frac{d c_m^-}{d
y} = - i \frac{\partial H}{\partial \bar{c}_m^-}, \quad m \in
\mathbb{Z}'_-
\end{equation}
and
\begin{equation}
\label{Ham-system-transformed-2} \frac{d c_m^+}{d y} = -
\frac{\partial H}{\partial \bar{c}_m^-}, \quad  \frac{d c_m^-}{d
y} = \frac{\partial H}{\partial \bar{c}_m^+}, \quad m \in
\mathbb{Z}'_+.
\end{equation}
Since the convolution sums on $\{ c_m^{\pm} \}_{m \in
\mathbb{Z}'}$ involve decaying weights as $|m| \to \infty$ and the
linear unbounded part on $\{ c_m\}_{m \in \mathbb{Z}'}$ is linear
in $|m|$ ($k_m^{\pm}$ and $\kappa_m^{\pm}$ grow linearly as $|m|
\to \infty$), the vector field of the system
(\ref{Ham-system-transformed-1})--(\ref{Ham-system-transformed-2})
has a modified domain $D'$ and range $X'$, which are given by
\begin{eqnarray}
D' = \left\{ ({\bf c}^+,{\bf c}^-,\bar{\bf c}^+,\bar{\bf c}^-) \in
l^2_{s'+1}(\mathbb{Z}',\mathbb{C}^4) \right\}, \qquad X' = \left\{
({\bf c}^+,{\bf c}^-,\bar{\bf c}^+,\bar{\bf c}^-) \in
l^2_{s'}(\mathbb{Z}',\mathbb{C}^4) \right\},
\label{new-domain-range}
\end{eqnarray}
for $s' = s - \frac{1}{4}$. (Here we have used the fact that the
convolution sum acts on $\{ \psi_m\}_{m \in \mathbb{Z}}$, where
$\psi_m$ is given by $c_m^{\pm}$ in the transformation
(\ref{normal-coordinates-1})--(\ref{normal-coordinates-2}).) If
the space $X$ is Banach algebra for $s > \frac{1}{2}$, the space
$X'$ is Banach algebra with respect to the decaying weights for
$s' > \frac{1}{4}$. The domain for reversible solutions become now
\begin{equation}
\label{domain-reversibility-modified} D_r' = \left\{ ({\bf
c}^+,{\bf c}^-,\bar{\bf c}^+,\bar{\bf c}^-) \in D' : \quad
c_m^{\pm}(-y) = \bar{c}_m^{\pm}(y), \; m \in \mathbb{Z}_-', \quad
c_m^{\pm}(-y) = \bar{c}_m^{\mp}(y), \; m \in \mathbb{Z}_+',
\right\}
\end{equation}
and the reversibility constraint becomes
\begin{equation}
\label{constraint-reversibility-modified} \Sigma_r' = \left\{
({\bf c}^+,{\bf c}^-,\bar{\bf c}^+,\bar{\bf c}^-) \in D' : \;\;
{\rm Im} c_m^{\pm} = 0, \; m \in \mathbb{Z}_-', \;\; {\rm Re}
c_m^+ = {\rm Re} c_m^-, \; {\rm Im} c_m^+ = - {\rm Im} c_m^-, \; m
\in \mathbb{Z}'_+ \right\}.
\end{equation}
In what follows, we are concerned with the reversible homoclinic
orbit of the Hamilton's equations of motion
(\ref{Ham-system-transformed-1})--(\ref{Ham-system-transformed-2})
in normal coordinates
(\ref{normal-coordinates-1})--(\ref{normal-coordinates-2}). We
will suppress the explicit dependence of $H$ from $(\bar{\bf
c}^+,\bar{\bf c}^-)$ for briefness of notations.

\begin{remark}
{\rm Since $k_m^{\pm} > 0$ for sufficiently large negative $m \in
\mathbb{Z}'_-$, the center manifold $E_c^+ \oplus E_c^-$ is
spanned by an infinite set of modes with positive and negative
energies. Therefore, we can not use the technique of \cite{GS01}
which relies on the fact that the quadratic part of the
Hamiltonian function is positive-definite for the non-bifurcating
modes of center manifold. We can however use the technique of
\cite{GS05} which relies on the separation of slow motion for the
modes, which correspond to the zero eigenvalue (modes $c_n^+$ and
$c_{-n}^-$), and the fast motion of the other modes, which
correspond to non-zero purely imaginary or complex eigenvalues
(modes $c_m^{\pm}$ for all other $m$). Moreover, we can simplify
the technique of \cite{GS05} by incorporating the Hamiltonian
structure
(\ref{Ham-system-transformed-1})--(\ref{Ham-system-transformed-2})
with the Hamiltonian function (\ref{Hamiltonian-transformed}). }
\end{remark}

\section{Normal form and persistence of homoclinic orbits}

We first show that the formal truncation of the Hamiltonian
function (\ref{Hamiltonian-transformed}) at the two bifurcating
modes $c_n^+$ and $c_{-n}^-$, which correspond to the double zero
eigenvalue of the linearized system at the zero solution for
$\epsilon = 0$, leads to the coupled-mode system (\ref{cme}) for
$\epsilon \neq 0$. Then, we derive an extended coupled-mode system
for the modes $c_n^+$ and $c_{-n}^-$ by using near-identity
transformations of the Hamiltonian function and prove persistence
of a reversible homoclinic orbit in the extended coupled-mode
system. We will assume from now on that $0 < c < n$.

Let us consider the subspace of the phase space of the Hamiltonian
system
(\ref{Ham-system-transformed-1})--(\ref{Ham-system-transformed-2}):
\begin{equation}
\label{invariant-subspace} S = \left\{ c_m^+ = 0, \; \forall m \in
\mathbb{Z}' \backslash \{n\}, \quad c_m^- = 0, \; \forall m \in
\mathbb{Z}' \backslash \{-n\} \right\}.
\end{equation}
Denote $Q = S^{\perp}$. If the Hamiltonian function $H$ is
formally constrained on the subspace $S$, the expression
(\ref{Hamiltonian-transformed}) takes the form
\begin{eqnarray*} H |_S = \epsilon \left[ \frac{\Omega
|c_n^+|^2}{n-c} + \frac{\Omega |c_{-n}^-|^2}{n+c} - \frac{w_{2n}
(\bar{c}_n^+ c_{-n}^- + c_{n}^+ \bar{c}_{-n}^-)}{\sqrt{n^2-c^2}} -
\frac{\sigma}{2} \left( \frac{|c_n^+|^4}{(n-c)^2} + \frac{4
|c_n^+|^2 |c_{-n}^-|^2}{n^2 - c^2} + \frac{|c_{-n}^-|^4}{(n+c)^2}
\right) \right].
\end{eqnarray*}
Since $(n,-n) \in \mathbb{Z}'_-$, we use the symplectic structure
(\ref{Ham-system-transformed-1}) to generate a system of
first-order ODEs for normal coordinates $c_n^+$ and $c_{-n}^-$. By
using the new independent variable $Y = \epsilon y$ and the new
dependent variables $a = \frac{c_n^+}{\sqrt{n-c}}$ and $b =
\frac{c_{-n}^-}{\sqrt{n+c}}$, we recover the ODE system
\begin{equation}
\label{cme-ode-system} \left\{ \begin{array}{ccc} i (n-c) a' +
\Omega a & = & w_{2n} b + \sigma (|a|^2 + 2 |b|^2) a, \\ - i (n+c)
b' + \Omega b & = & w_{2n} a + \sigma (2 |a|^2 + |b|^2) b,
\end{array} \right.
\end{equation}
where the derivatives are taken with respect to $Y = \epsilon y$.
The system (\ref{cme-ode-system}) is nothing but the coupled-mode
system (\ref{cme}) after the separation of the traveling variable
$Y = X - cT$ and the stationary variable $T$ in the transformation
$$
a(X,T) \mapsto a(Y) e^{-i \Omega T}, \quad b(X,T) \mapsto b(Y)
e^{-i \Omega T},  \quad Y = X-cT,
$$
with the correspondence $Y = \sqrt{n^2-c^2} \xi$ and $\Omega =
\frac{\mu n (1-c^2)}{\sqrt{n^2-c^2}}$. The system
(\ref{cme-ode-system}) has a localized solution (a homoclinic
orbit) (see Section 1) for $w_{2n} \neq 0$ and $|\Omega| <
\Omega_0$, where $\Omega_0 = |w_{2n}| \frac{\sqrt{n^2 - c^2}}{n}$.

The ODE system (\ref{cme-ode-system}) is invariant with respect to
translation $a(Y) \to a(Y-Y_0)$, $b(Y) \to b(Y-Y_0)$ for all $Y_0
\in \mathbb{R}$ and gauge transformation $a(Y) \to e^{i \alpha}
a(Y-Y_0)$, $b(Y) \to e^{i \alpha} b(Y-Y_0)$ for all $\alpha \in
\mathbb{R}$. Therefore, any solution of the system is continued
with a two-parameter group of symmetry transformations. However,
these parameters are set uniquely in the reversible homoclinic
orbit of Definition \ref{definition-orbit}. In addition, we note
that although the ODE system (\ref{cme-ode-system}) is formulated
in the four-dimensional phase space, it has two conserved
quantities on $Y \in \mathbb{R}$ related to the translational and
gauge symmetries. Indeed, the Hamiltonian $H |_S$ and the
quadratic function
\begin{equation}
Q = |c_n^+|^2 - |c_{-n}^-|^2 = (n-c)|a|^2 - (n+c)|b|^2
\end{equation}
are constants on $Y \in \mathbb{R}$. As a result, localized
solutions of the system (\ref{cme-ode-system}) are defined on a
subspace which obeys a planar Hamiltonian system. This planar
Hamiltonian system is given by the second equation of the system
(\ref{varphi-phi}) in variables $(\phi,\bar{\phi})$. We recall
that reversible homoclinic orbits of planar Hamiltonian systems
are structurally stable with respect to parameter continuations.

To incorporate the ideas of integrability of the Hamiltonian
coupled-mode system on the subspace $S$ and persistence of the
reversible homoclinic orbit in the planar Hamiltonian system, we
extend the coupled-mode system by using near-identity
transformations and the normal form theory.

\begin{lemma}
Let $0 < c < n$, such that $\frac{n^2+c^2}{2c} \notin
\mathbb{Z}'$. For each $N \in \mathbb{N}$ and sufficiently small
$\epsilon$, there is a near-identity, analytic, symplectic change
of coordinates in a neighborhood of the origin in $X'$ in
(\ref{new-domain-range}), such that the Hamiltonian function $H$
in (\ref{Hamiltonian-transformed}) transforms to the normal form
up to the order of ${\rm O}(\epsilon^{N+1})$,
\begin{eqnarray}
\nonumber H & = & \sum_{m \in \mathbb{Z}'_-} \left( k_m^+
|c_m^+|^2 - k_m^- |c_m^-|^2 \right) + \sum_{m \in \mathbb{Z}'_+}
\left( \kappa_m^- c_m^- \bar{c}_m^+ - \kappa_m^+
c_m^+ \bar{c}_m^- \right) \\
\label{Hamiltonian-normal-form} & \phantom{t} & \phantom{text} +
\epsilon H_S(c_n^+,c_{-n}^-) + \epsilon H_Q(c_n^+,c_{-n}^-,{\bf
c}^+,{\bf c}^-) + \epsilon^{N+1} H_R(c_n^+,c_{-n}^-,{\bf c}^+,{\bf
c}^-),
\end{eqnarray}
where $H_S$ is a polynomial of the degree $2N+2$ in
$(c_n^+,c_{-n}^-)$ on $S$, $H_Q$ is a polynomial of the degree
$2N$ in $(c_n^+,c_{-n}^-)$ on $S$ and of the degree $4$ in $({\bf
c}^+,{\bf c}^-)$ on $Q$ with no linear terms with respect to
$({\bf c}^+,{\bf c}^-)$ on $Q$, and $H_R$ is a polynomial of the
degree $8 N + 4$ in $(c_n^+,c_{-n}^-)$ on $S$ and of the degree 4
in $({\bf c}^+,{\bf c}^-)$ on $Q$. All components $H_S$, $H_Q$ and
$H_R$ depend on $\epsilon$, such that $H_S$ and $H_Q$ are
polynomials in $\epsilon$ of the degree $N-1$ and $H_R$ is a
polynomial in $\epsilon$ of the degree $3N-1$. The reversibility
(\ref{transformation-reversibility}) and gauge
(\ref{transformation-gauge}) transformations are preserved by the
change of the variables. \label{lemma-normal-form-transformation}
\end{lemma}

\begin{proof}
The existence of a near-identity symplectic transformation that
maps $H$ to the form (\ref{Hamiltonian-normal-form}) follows from
the fact that the non-resonance conditions $l \kappa_0 -
\kappa_m^{\pm} \neq 0$ are satisfied for any $l \in \mathbb{Z}$
and $m \in \mathbb{Z}'$ since $\kappa_0 = 0$ and all eigenvalues
are semi-simple for $\epsilon = 0$. See \cite{GS01} for an
iterative sequence of symplectic transformations. The
transformation is analytic in a local neighborhood of the origin
in $X'$ as the vector field of the Hamiltonian system
(\ref{Ham-system-transformed-1})--(\ref{Ham-system-transformed-2})
is analytic (given by a cubic polynomial). The reversibility
(\ref{transformation-reversibility}) and gauge
(\ref{transformation-gauge}) symmetries are preserved by the
symplectic change of variables \cite{GS01}. The count of the
degree of polynomials $H_S$, $H_Q$ and $H_R$ follows from the fact
that the vector field of the Hamiltonian system
(\ref{Ham-system-transformed-1})--(\ref{Ham-system-transformed-2})
contains only linear and cubic terms in normal coordinates, while
the near-identity transformation of $({\bf c}^+,{\bf c}^-)$ on $Q$
up to the order of ${\rm O}(\epsilon^{N+1})$ involves a polynomial
in $\epsilon$ of the degree $N$ and a polynomial in
$(c_n^+,c_{-n}^-)$ on $S$ of the degree $2N+1$.
\end{proof}

\begin{example}
{\rm For $N = 1$, the near-identity transformation for $c_m^+$, $m
\in \mathbb{Z}'_- \backslash \{ n\}$ takes the explicit form
$$
\tilde{c}_m^+ = c_m^+ - \frac{\epsilon}{k_m^+} \left[ w_{m,n}
c_n^+ + w_{m,-n} c_{-n}^- + \sigma \left( g_{n,n,n} (c_n^+)^2
\bar{c}_{-n}^- \delta_{m,3n} + g_{-n,-n,-n} |c_n^+|^2 c_{-n}^-
\delta_{m,-n} \right) \right] + {\rm O}(\epsilon^2),
$$
where $\tilde{c}_m^+$ is a new variable. Similar explicit formulas
can be obtained for $c_m^-$, $m \in \mathbb{Z}_-' \backslash
\{-n\}$ and for $c_m^{\pm}$, $m \in \mathbb{Z}'_+$. We note that
the tilde signs are omitted for new variables in the expression
(\ref{Hamiltonian-normal-form}).}
\label{example-normal-transformation}
\end{example}

\begin{remark}
\label{remark-extended-cme} {\rm For each $N \in \mathbb{N}$, the
subspace $S$ defined by (\ref{invariant-subspace}) is an invariant
subspace of the Hamiltonian system
(\ref{Ham-system-transformed-1})--(\ref{Ham-system-transformed-2})
with the Hamiltonian function (\ref{Hamiltonian-normal-form})
truncated at $H_R \equiv 0$. The dynamics on $S$ is given by the
four-dimensional Hamiltonian system
\begin{equation}
\frac{d c_n^+}{dY} = i \frac{\partial H_S}{\partial \bar{c}_n^+},
\qquad \frac{d c_{-n}^-}{dY} = -i \frac{\partial H_S}{\partial
\bar{c}_{-n}^+}, \label{cme-normal-form}
\end{equation}
where $Y = \epsilon y$. If $N = 1$, the system
(\ref{cme-normal-form}) transforms to the coupled-mode system
(\ref{cme-ode-system}) in variables $a = \frac{c_n^+}{\sqrt{n-c}}$
and $b = \frac{c_{-n}^-}{\sqrt{n+c}}$. If $N > 1$, this system is
referred to as the {\em extended} coupled-mode system. }
\end{remark}

\begin{lemma}
Let $w_{2n} \neq 0$ for a given $n \in \mathbb{N}$. For each $N
\in \mathbb{N}$ and sufficiently small $\epsilon$, there exists a
reversible homoclinic orbit of the system (\ref{cme-normal-form})
for $|\Omega| < \Omega_0 = |w_{2n}| \frac{\sqrt{n^2 - c^2}}{n}$.
Moreover, the solution for the homoclinic orbit satisfies the
global bound
\begin{equation}
|c_n^+(y)| \leq C_+ e^{-\epsilon \gamma |y|}, \quad  |c_{-n}^-(y)|
\leq C_- e^{-\epsilon \gamma |y|}, \quad \forall y \in \mathbb{R},
\label{bound-homoclinic-orbit}
\end{equation}
for some $\epsilon$-independent constants $\gamma
> 0$ and $C_{\pm} > 0$.
\label{lemma-homoclinic-orbit}
\end{lemma}

\begin{proof}
Due to the gauge-invariance of the polynomial Hamiltonian function
$H_S = H_S(c_n^+,c_{-n}^-)$, it must satisfy the partial
differential equation \cite{ChPel}:
\begin{equation}
\frac{d}{d \alpha} H_S(e^{i \alpha} c_n^+, e^{i \alpha} c_{-n}^-)
\biggr|_{\alpha = 0} \backsimeq c_n^+ \frac{\partial H_S}{\partial
c_n^+} - \bar{c}_n^+ \frac{\partial H_S}{\partial \bar{c}_n^+} +
c_{-n}^- \frac{\partial H_S}{\partial c_{-n}^-} - \bar{c}_{-n}^-
\frac{\partial H_S}{\partial c_{-n}^-} = 0. \label{gauge-relation}
\end{equation}
It follows from the system (\ref{cme-normal-form}) and the
relation (\ref{gauge-relation}) that $Q = |c_n^+|^2 -
|c_{-n}^-|^2$ is constant in $Y \in \mathbb{R}$ \cite{ChPel}. If
localized solutions exist, then $Q = 0$. Let us represent the
solutions in the general form
\begin{equation}
\label{extended-normal-coordinates} c_n^+ = \sqrt{\rho + q} e^{i
\varphi + i \theta}, \quad c_{-n}^- = \sqrt{\rho - q} e^{i \varphi
- i \theta},
\end{equation}
where $(\rho,q,\theta,\varphi)$ are new real-valued variables.
Using the chain rule for $H_S =
\tilde{H}_S(\rho,q,\theta,\varphi)$, we find that $\varphi$ is a
cyclic variable of the Hamiltonian function $\tilde{H}_S$ and $q$
is constant due to the gauge invariance (\ref{gauge-relation}).
Setting $q = 0$, we find that $(\rho,\theta)$ satisfy a planar
Hamiltonian system, while $\varphi$ is found from a linear
inhomogeneous equation:
\begin{equation}
\label{extended-cme} 2 \frac{d \rho}{d Y} = - \frac{\partial
\tilde{H}_S}{\partial \theta} |_{q = 0}, \qquad 2 \frac{d
\theta}{d Y} = \frac{\partial \tilde{H}_S}{\partial \rho} |_{q =
0}, \qquad 2 \frac{d \varphi}{d Y} = \frac{\partial
\tilde{H}_S}{\partial q} |_{q = 0},
\end{equation}
where $\tilde{H}_S |_{q = 0}$ is a function of $(\rho,\theta)$. If
$\epsilon = 0$, the system (\ref{extended-cme}) reduces to the ODE
system (\ref{varphi-phi}) rewritten in new coordinates and for
$\phi = \sqrt{\rho} e^{i \theta}$. The vector field of the
extended coupled-mode system (\ref{cme-normal-form}) is given by
polynomials in $c_n^+$ and $c_{-n}^-$ of the degree $2N + 1$ and
in $\epsilon$ of the degree $N-1$. Recall that the coupled-mode
system (\ref{cme-ode-system}) admits a reversible homoclinic orbit
for $w_{2n} \neq 0$ and $|\Omega| < \Omega_0$, where $\Omega_0 =
|w_{2n}| \frac{\sqrt{n^2-c^2}}{n}$. Since a reversible homoclinic
orbit is structurally stable in a planar Hamiltonian system with
an analytic vector field, the homoclinic orbit persists in the
extended coupled-mode system (\ref{cme-normal-form}) for
sufficiently small $\epsilon$.
\end{proof}

\section{Construction of local center and center-stable manifolds}

We study solutions of the Hamiltonian system of equations
(\ref{Ham-system-transformed-1})--(\ref{Ham-system-transformed-2})
after the normal-form transformation of Lemma
\ref{lemma-normal-form-transformation}. We construct a local
solution on $y \in [0,L/\epsilon^{N+1}]$ for some
$\epsilon$-independent constant $L > 0$, which is close to the
homoclinic orbit of Lemma \ref{lemma-homoclinic-orbit} by the
distance $C \epsilon^N$ for some $\epsilon$-independent constant
$C > 0$ in vector norm on $X'$. This solution represents an
infinite-dimensional local center--stable manifold and it is
spanned by the small bounded oscillatory and small exponentially
decaying solutions near the exponentially decaying homoclinic
solution with the decay bound (\ref{bound-homoclinic-orbit}).
Parameters of the local center--stable manifold are chosen to
ensure that the manifold intersects at $y = 0$ with the symmetric
section $\Sigma_r'$. This construction completes the proof of the
main Theorem \ref{theorem-main}.

By using Lemma \ref{lemma-normal-form-transformation} and the
explicit representation of the Hamiltonian function
(\ref{Hamiltonian-normal-form}), we rewrite the Hamiltonian system
of equations in the separated form
\begin{eqnarray}
\label{semilinear-1} \frac{d c_n^+}{d y} & = & \epsilon \left[
\mu^+_{\epsilon} c_n^+ + \nu^+_{\epsilon} c_{-n}^- +
F_S^+(c_n^+,c_{-n}^-) + F_Q^+(c_n^+,c_{-n}^-,{\bf c}) \right] +
\epsilon^{N+1} F_R^+(c_n^+,c_n^-,{\bf c}), \\
\label{semilinear-2} \frac{d c_{-n}^-}{d y} & = & \epsilon \left[
\mu^-_{\epsilon} c_n^+ + \nu^-_{\epsilon} c_{-n}^- +
F_S^-(c_n^+,c_{-n}^-) + F_Q^-(c_n^+,c_{-n}^-,{\bf c}) \right] +
\epsilon^{N+1} F_R^-(c_n^+,c_n^-,{\bf c}), \\
\label{semilinear-3} \frac{d {\bf c}}{d y} & = &
\Lambda_{\epsilon} {\bf c} + \epsilon {\bf
F}_Q(c_n^+,c_{-n}^-,{\bf c}) + \epsilon^{N+1} {\bf
F}_R(c_n^+,c_n^-,{\bf c}),
\end{eqnarray}
where ${\bf c}$ denotes all components of $({\bf c}^+,{\bf
c}^-,\bar{\bf c}^+,\bar{\bf c}^-)$ in $Q = S^{\perp}$,
($\mu^{\pm}_{\epsilon}$,$\nu_{\epsilon}^{\pm}$,$\Lambda_{\epsilon}$)
denote the coefficient matrix for the linear part of the system
and $(F_S^{\pm},F_Q^{\pm},F_R^{\pm},{\bf F}_Q,{\bf F}_R)$ denote
the nonlinear (polynomial) part of the system. For briefness of
notations, we do not rewrite the subsystem
(\ref{semilinear-1})--(\ref{semilinear-2}) for variables
$(\bar{c}_n^+,\bar{c}_{-n}^-)$ and we do not write dependence of
the nonlinear functions from these variables. The variables ${\bf
c}$ are equivalent to tilde-variables in Example
\ref{example-normal-transformation} after the near-identity
transformations, but the tilde-notations are dropped for
simplicity of notations.

\begin{lemma}
\label{lemma-perturbation-theory} Let $W(x)$ satisfy Assumption
\ref{assumption-potential} and $w_{2n} \neq 0$ for a given $n \in
\mathbb{N}$. For sufficiently small $\epsilon$, the linearized
system (\ref{semilinear-1})--(\ref{semilinear-3}) at the zero
solution for $|\Omega| < \Omega_0 = |w_{2n}| \frac{\sqrt{n^2 -
c^2}}{n}$ and $0 < c < n$ is topologically equivalent to the one
for $\epsilon = 0$, except that the double zero eigenvalue of the
subsystem (\ref{semilinear-1})--(\ref{semilinear-2}) splits into a
pair of complex eigenvalues to the left and right half-planes.
\end{lemma}

\begin{proof}
Since all non-zero eigenvalues of the linearized Hamiltonian
system (\ref{semilinear-3}) at the zero solution are semi-simple
at $\epsilon = 0$, they are structurally stable in the
perturbation theory for sufficiently small $\epsilon \neq 0$. The
matrix operator $\Lambda_{\epsilon}$ is a polynomial in $\epsilon$
and $\Lambda_0$ is a diagonal unbounded matrix operator which
consists of $(i k_m^+,-i k_m^-,-ik_m^+,ik_m^-)$ for $m \in
\mathbb{Z}'_-$ and of
$(\kappa_m^+,\kappa_m^-,\bar{\kappa}_m^+,\bar{\kappa}_m^-)$ for $m
\in \mathbb{Z}'_+$. The matrix operator with elements $\epsilon
w_{m,m'}$ represents a small perturbation to $\Lambda_0$ if the
vector of Fourier coefficients ${\bf W}$ is in $l^1(\mathbb{Z})
\subset l^2_s(\mathbb{Z})$ for $s
> \frac{1}{2}$ according to Assumption \ref{assumption-potential}.

\noindent The coefficients $\mu_{\epsilon}^{\pm}$ and
$\nu_{\epsilon}^{\pm}$ are polynomials in $\epsilon$ and
$$
\mu^+_0 = \frac{i \Omega}{n - c}, \quad \mu^-_0 = \frac{i
w_{2n}}{\sqrt{n^2 - c^2}}, \quad \nu^+_0 = \frac{-i
w_{2n}}{\sqrt{n^2 - c^2}}, \quad \nu^-_0 = \frac{-i \Omega}{n +
c}.
$$
At $\epsilon = 0$, the linearized subsystem
(\ref{semilinear-1})--(\ref{semilinear-2}) corresponds to the
linearized coupled-mode system (\ref{cme-ode-system}). Its
characteristic equation is given by
$$
(n^2-c^2) \kappa^2 - 2 i \epsilon c \Omega \kappa + \epsilon^2
(\Omega^2 - w_{2n}^2) = 0,
$$
with two roots
$$
\kappa = \kappa_{\pm} = \epsilon \frac{i \Omega c \pm
\sqrt{(n^2-c^2) w_{2n}^2 - n^2 \Omega^2}}{n^2 - c^2}.
$$
The two roots have ${\rm Re} \kappa_{\pm} \gtrless 0$ if $|\Omega|
< \Omega_0$, where $\Omega_0 = |w_{2n}| \frac{\sqrt{n^2 -
c^2}}{n}$. Under the same assumption on ${\bf W}$, perturbation
terms in $\mu_{\epsilon}^{\pm}$ and $\nu_{\epsilon}^{\pm}$ are
small compared to the leading-order terms $\mu_0^{\pm}$ and
$\nu_0^{\pm}$, such that the pair persists in the left half-plane
and right half-plane of the $\kappa$-plane.
\end{proof}

\begin{corollary}
\label{corollary-decomposition} For sufficiently small $\epsilon$,
a local neighborhood of the zero point in the phase space $X'$ can
be decomposed into the subspaces determined by the spectrum of the
linearized system at the zero solution
\begin{equation}
\label{decomposition-final} X' = X_h \oplus X_c \oplus X_u \oplus
X_s,
\end{equation}
where $X_s$ and $X_u$ are associated to the subspaces $E^s$ and
$E^u$ of Lemma \ref{lemma-roots}, while $X_h$ and $X_c$ are
associated to the subspaces $E^{c^+} \oplus E^{c^-}$ on $S$ and
$Q$ respectively.
\end{corollary}

By Remark \ref{remark-extended-cme}, the truncated system
(\ref{semilinear-1})--(\ref{semilinear-3}) with $F_R^{\pm} \equiv
0$ and ${\bf F}_R \equiv {\bf 0}$ admits an invariant reduction on
$S$. By Lemma \ref{lemma-homoclinic-orbit}, the extended
coupled-mode system (\ref{cme-normal-form}) on $S$ has a
reversible homoclinic orbit which satisfies the decay bound
(\ref{bound-homoclinic-orbit}). This construction enables us to
represent the solution of the subsystem
(\ref{semilinear-1})--(\ref{semilinear-2}) and its complex
conjugate in the form
$[c_n^+,c_{-n}^-,\bar{c}_n^+,\bar{c}_{-n}^-]^T = {\bf
c}_0(\epsilon y) + {\bf c}_h(y)$, where ${\bf c}_0(\epsilon y)$ is
the homoclinic orbit of Lemma \ref{lemma-homoclinic-orbit} and
${\bf c}_h(y)$ is a perturbation term. By using the decomposition,
we rewrite the system (\ref{semilinear-1})--(\ref{semilinear-3})
in the equivalent form
\begin{eqnarray}
\label{semilinear-4} \frac{d {\bf c}_h}{d y} & = & \epsilon
\Lambda_h({\bf c}_0) {\bf c}_h + \epsilon {\bf G}_Q({\bf
c}_0)({\bf c}_h,{\bf c}) +
\epsilon^{N+1} {\bf G}_R({\bf c}_0 + {\bf c}_h,{\bf c}), \\
\label{semilinear-5} \frac{d {\bf c}}{d y} & = &
\Lambda_{\epsilon} {\bf c} + \epsilon {\bf F}_Q({\bf c}_0 + {\bf
c}_h,{\bf c}) + \epsilon^{N+1} {\bf F}_R({\bf c}_0 + {\bf
c}_h,{\bf c}),
\end{eqnarray}
where $\Lambda_h({\bf c}_0)$ is a $4$-by-$4$ linearization matrix
of the extended coupled-mode system (\ref{cme-normal-form}) and
its complex conjugate at the solution ${\bf c}_0(\epsilon y)$ and
$({\bf G}_Q,{\bf G}_R)$ denote the nonlinear part of the subsystem
(\ref{semilinear-1})--(\ref{semilinear-2}) and its conjugate. We
note that the function ${\bf G}_Q$ combines nonlinear terms in
${\bf c}_h$ from the functions $F_S^{\pm}$ and the nonlinear terms
in ${\bf c}$ from the functions $F_Q^{\pm}$. For simplicity of
notations, we say that ${\bf c}_h \in X_h$ and ${\bf c} \in
X_h^{\perp}$ in the decomposition $X' = X_h \oplus X_h^{\perp}$.

\begin{example}
{\rm For $N = 1$, the linearization matrix $\Lambda_h({\bf c}_0)$
takes the explicit form
$$
i \left[ \begin{array}{cccc} \frac{\Omega}{n-c} - \frac{2 \sigma
|c_n^+|^2}{(n-c)^2} - \frac{2 \sigma |c_{-n}^-|^2}{n^2-c^2} &
-\frac{w_{2n}}{\sqrt{n^2-c^2}} - \frac{2 \sigma c_n^+
\bar{c}_{-n}^-}{n^2-c^2} & -\frac{\sigma c_n^{+2}}{(n-c)^2} &
-\frac{2 \sigma c_n^+ c_{-n}^-}{n^2-c^2} \\
\frac{w_{2n}}{\sqrt{n^2-c^2}} + \frac{2 \sigma \bar{c}_n^+
c_{-n}^-}{n^2-c^2} & - \frac{\Omega}{n+c} + \frac{2 \sigma
|c_n^+|^2}{n^2-c^2} + \frac{2 \sigma |c_{-n}^-|^2}{(n+c)^2} &
\frac{2 \sigma c_n^+ c_{-n}^-}{n^2-c^2} & \frac{\sigma c_{-n}^{-2}}{(n+c)^2} \\
\frac{\sigma \bar{c}_n^{+2}}{(n-c)^2} & \frac{2 \sigma \bar{c}_n^+
\bar{c}_{-n}^-}{n^2-c^2} & -\frac{\Omega}{n-c} + \frac{2 \sigma
|c_n^+|^2}{(n-c)^2} + \frac{2 \sigma |c_{-n}^-|^2}{n^2-c^2} &
\frac{w_{2n}}{\sqrt{n^2-c^2}} +
\frac{2 \sigma \bar{c}_n^+ c_{-n}^-}{n^2-c^2} \\
-\frac{2 \sigma \bar{c}_n^+ \bar{c}_{-n}^-}{n^2-c^2} & -
\frac{\sigma \bar{c}_{-n}^{-2}}{(n+c)^2} &
-\frac{w_{2n}}{\sqrt{n^2-c^2}} - \frac{2 \sigma c_n^+
\bar{c}_{-n}^-}{n^2-c^2} & \frac{\Omega}{n+c} - \frac{2 \sigma
|c_n^+|^2}{n^2-c^2} - \frac{2 \sigma |c_{-n}^-|^2}{(n+c)^2}
\end{array} \right],
$$
where components of
$[c_n^+,c_{-n}^-,\bar{c}_n^+,\bar{c}_{-n}^-]^T$ are evaluated at
the solution ${\bf c}_0$ of the extended coupled-mode equation
(\ref{cme-normal-form}). As $|y| \to \infty$, the matrix
$\Lambda_h({\bf c}_0)$ converges to the form
$$
\Lambda_h({\bf 0}) = i \left[ \begin{array}{cccc}
\frac{\Omega}{n-c} &
-\frac{w_{2n}}{\sqrt{n^2-c^2}}  & 0 & 0 \\
\frac{w_{2n}}{\sqrt{n^2-c^2}} & - \frac{\Omega}{n+c} & 0 & 0 \\
0 & 0 & -\frac{\Omega}{n-c} & \frac{w_{2n}}{\sqrt{n^2-c^2}} \\
0 & 0 & -\frac{w_{2n}}{\sqrt{n^2-c^2}} & \frac{\Omega}{n+c}
\end{array} \right],
$$
such that $\| \Lambda_h({\bf c}_0) - \Lambda_h({\bf 0}) \|_{X_h
\mapsto X_h} \leq C e^{-\epsilon \gamma |y|}$ for some $C > 0$ and
$\gamma > 0$ according to the decay bound
(\ref{bound-homoclinic-orbit}). }
\label{example-linearized-operator}
\end{example}

\begin{lemma}
Let $w_{2n} \neq 0$ for a given $n \in \mathbb{N}$, $|\Omega| <
\Omega_0 = |w_{2n}| \frac{\sqrt{n^2 - c^2}}{n}$, and $0 < c < n$.
Consider the linear inhomogeneous equation
\begin{equation}
\label{inhomogeneous-equation-h} \frac{d {\bf c}_h}{d y} -
\epsilon \Lambda_h({\bf c}_0) {\bf c}_h = {\bf F}_h(y),
\end{equation}
where ${\bf F}_h \in C_b^0(\mathbb{R})$. The homogeneous equation
has a two-dimensional stable manifold spanned by the two
fundamental solutions
\begin{equation}
\label{kernel} {\bf s}_1 = {\bf c}_0'(y), \quad {\bf s}_2 =
\sigma_1 {\bf c}_0(y),
\end{equation}
where $\sigma_1$ is a diagonal matrix of $(1,1,-1,-1)$. If
components of ${\bf F}_h(y)$ satisfies the constraints $\forall y
\in \mathbb{R}$
\begin{equation}
\label{constraints-kernel} ({\bf F}_h)_1(y) = (\bar{\bf
F}_h)_1(-y), \;\; ({\bf F}_h)_2(y) = (\bar{\bf F}_h)_2(-y), \;\;
({\bf F}_h)_3(y) = -(\bar{\bf F}_h)_1(y), \;\; ({\bf F}_h)_4(y) =
-(\bar{\bf F}_h)_4(y),
\end{equation}
then there exists a two-parameter family of solutions ${\bf c}_h
\in C_b^0(\mathbb{R})$ in the form ${\bf c}_h = \alpha_1 {\bf
s}_1(y) + \alpha_2 {\bf s}_2(y) + \tilde{\bf c}_h(y)$, where
$(\alpha_1,\alpha_2)$ are parameters and $\tilde{\bf c}_h(y)$ is a
particular solution of the inhomogeneous equation
(\ref{inhomogeneous-equation-h}), such that $\| \tilde{\bf c}_h
\|_{C_b^0(\mathbb{R})} \leq \frac{C}{\epsilon} \| {\bf F}_h
\|_{C^0_b(\mathbb{R})}$ for an $\epsilon$-independent constant
$C$. \label{lemma-weak-stable-unstable-manifolds}
\end{lemma}

\begin{proof}
The existence of the two-dimensional kernel (\ref{kernel}) follows
from symmetries of the extended coupled-mode system
(\ref{cme-normal-form}) with respect to translation and gauge
transformation. Since the subspace $X_h$ associated with
$\Lambda_h({\bf 0})$ is invariant under $\Lambda_h({\bf c}_0)$,
the kernel is exactly two-dimensional and the other two
fundamental solutions of the homogeneous equation are
exponentially growing. The adjoint homogeneous equation $\frac{d
{\bf u}}{d y} = - \epsilon \Lambda_h^*({\bf c}_0) {\bf u}$ has
also a two-dimensional stable manifold spanned by the two
fundamental solutions
\begin{equation}
\label{kernel-adjoint} {\bf s}^*_1 = \sigma_2 {\bf c}_0'(y), \quad
{\bf s}^*_2 = \sigma_3 {\bf c}_0(y),
\end{equation}
where $\sigma_2$ and $\sigma_3$ are diagonal matrices of
$(1,-1,-1,1)$ and $(1,-1,1,-1)$ respectively. Unless the vector
function ${\bf F}_h$ is orthogonal to $\{ {\bf s}_1^*,{\bf s}_2^*
\}$, a solution of the linear inhomogeneous equation
(\ref{inhomogeneous-equation-h}) grows exponentially as $|y| \to
\infty$. However, $({\bf s}_1^*,{\bf F}_h) = ({\bf s}_2^*,{\bf
F}_h) = 0$ if the constraints (\ref{constraints-kernel}) are
satisfied. By the Fredholm theory, ${\bf F}_h$ is in the range of
the linear unbounded operator $\frac{d}{dy} - \epsilon
\Lambda_h({\bf c}_0)$, such that there exists a solution
$\tilde{\bf c}_h \in C_b^0(\mathbb{R})$ of the inhomogeneous
equation (\ref{inhomogeneous-equation-h}) such that $\| \tilde{\bf
c}_h \|_{C_b^0(\mathbb{R})} \leq \frac{C}{\epsilon} \| {\bf F}_h
\|_{C^0_b(\mathbb{R})}$ for an $\epsilon$-independent constant
$C$. A general solution of the inhomogeneous problem has the form
${\bf c}_h = \alpha_1 {\bf s}_1(y) + \alpha_2 {\bf s}_2(y) +
\tilde{\bf c}_h(y)$, where $(\alpha_1,\alpha_2)$ are parameters.
\end{proof}

\begin{lemma}
\label{lemma-nonlinear-vector-field} The nonlinear part of the
vector field of the system
(\ref{semilinear-4})--(\ref{semilinear-5}) is bounded in a local
neighborhood of the zero point in $X' = X_h \oplus X_h^{\perp}$ by
\begin{eqnarray}
\label{estimates-1} \| {\bf G}_R \|_{X_h} \leq N_R \left( \| {\bf
c}_0 + {\bf c}_h \|_{X_h} + \| {\bf c} \|_{X_h^{\perp}} \right),
\quad \| {\bf F}_R \|_{X^{\perp}_h} \leq M_R \left(  \| {\bf c}_0
+ {\bf c}_h \|_{X_h} + \| {\bf c} \|_{X_h^{\perp}}  \right), \\
\label{estimates-2} \| {\bf G}_Q \|_{X_h} \leq N_Q \left( \| {\bf
c}_h \|^2_{X_h} + \| {\bf c} \|^2_{X_h^{\perp}} \right), \quad \|
{\bf F}_Q \|_{X'} \leq M_Q  \left(  \| {\bf c}_0 + {\bf c}_h
\|_{X_h} + \| {\bf c} \|_{X_h^{\perp}}  \right) \| {\bf c}
\|_{X_h^{\perp}},
\end{eqnarray}
for some $N_R,M_R,N_Q,M_Q > 0$.
\end{lemma}

\begin{proof}
The system (\ref{semilinear-4})--(\ref{semilinear-5}) is
semi-linear with polynomial vector field for a finite $N \in
\mathbb{N}$ defined on the domain $D' \subset X'$, where $D'$ and
$X'$ are given in (\ref{new-domain-range}). If $X$ is the Banach
algebra for $s > \frac{1}{2}$, then $X'$ is Banach algebra for $s'
> \frac{1}{4}$. The derivation of estimates (\ref{estimates-1})--(\ref{estimates-2})
follows similarly to Lemmas 2 and 3 in \cite{GS01}. The
characterization of ${\bf G}_Q$ and ${\bf F}_Q$ is based on the
fact that the Hamiltonian function $H_Q$ is quadratic with respect
to ${\bf c}$ by Lemma \ref{lemma-normal-form-transformation}.
\end{proof}

\begin{remark}
{\rm Using Corollary \ref{corollary-decomposition}, we denote
$\Lambda_c = \Lambda_{\epsilon} |_{X_c}$, $\Lambda_u =
\Lambda_{\epsilon} |_{X_u}$, $\Lambda_s = \Lambda_{\epsilon}
|_{X_s}$ for a block-diagonal decomposition of
$\Lambda_{\epsilon}$ on the invariant subspaces $X_c$, $X_u$ and
$X_s$ respectively. We also denote the coordinates of the
decomposition by ${\bf c} = [{\bf c}_c,{\bf c}_u,{\bf c}_s]$ and
the projection operators by $P_c$, $P_u$, $P_s$ respectively. }
\end{remark}

\begin{theorem}{\bf (Local center-stable manifold)}
\label{theorem-reduction} Let ${\bf a} \in X_c$, ${\bf b} \in X_s$
and $(\alpha_1,\alpha_2) \in \mathbb{C}^2$ be small such that
\begin{equation}
\label{bound-initial} \|{\bf a} \|_{X_c} \leq C_a \epsilon^N,
\quad \| {\bf b} \|_{X_s} \leq C_b \epsilon^N, \quad |\alpha_1| +
|\alpha_2| \leq C_{\alpha} \epsilon^N.
\end{equation}
for some $\epsilon$-independent constants $C_a,C_b,C_{\alpha} >
0$. Under the conditions of Lemma \ref{lemma-perturbation-theory},
there exists a family ${\bf c}_h(y;{\bf a},{\bf
b},\alpha_1,\alpha_2)$ and ${\bf c}(y;{\bf a},{\bf
b},\alpha_1,\alpha_2)$ of local solutions of the system
(\ref{semilinear-4})--(\ref{semilinear-5}) such that ${\bf c}_c(0)
= {\bf a}$, ${\bf c}_s = e^{y \Lambda_s} {\bf b} + \tilde{\bf
c}_s(y)$, ${\bf c}_h = \alpha_1 {\bf s}_1(y) + \alpha_2 {\bf
s}_2(y) + \tilde{\bf c}_h(y)$ with uniquely defined $\tilde{\bf
c}_s(y)$ and $\tilde{\bf c}_h(y)$, and the local solutions satisfy
the bound
\begin{equation}
\label{bound-center} \sup_{y \in [0,L/\epsilon^{N+1}]} \| {\bf
c}_h(y) \|_{X_h} \leq C_h \epsilon^N, \qquad \sup_{y \in
[0,L/\epsilon^{N+1}]} \| {\bf c}(y) \|_{X_h^{\perp}} \leq C
\epsilon^N,
\end{equation}
for some $\epsilon$-independent constants $L > 0$ and $C_h,C > 0$.
\end{theorem}

\begin{proof}
We modify the system (\ref{semilinear-4})--(\ref{semilinear-5}) by
the following trick. We multiply the nonlinear vector field of the
subsystem (\ref{semilinear-5}) by the cut-off function
$\chi_{[0,y_0]}(y)$, such that
\begin{equation}
\label{semilinear-6} \frac{d {\bf c}}{d y} = \Lambda_{\epsilon}
{\bf c} + \epsilon \chi_{[0,y_0]}(y) {\bf F}_Q({\bf c}_0 + {\bf
c}_h,{\bf c}) + \epsilon^{N+1} \chi_{[0,y_0]}(y) {\bf F}_R({\bf
c}_0 + {\bf c}_h,{\bf c}),
\end{equation}
where $\chi_{[0,y_0]}(y) = 1$ for $y \in [0,y_0] \subset
\mathbb{R}$ for some $y_0 > 0$ and $\chi_{[0,y_0]}(y) = 0$
otherwise. Similarly, we multiply the nonlinear vector field of
the subsystem (\ref{semilinear-4}) by the cut-off function
$\chi_{[0,y_0]}(y)$ and add symmetrically reflected vector field
multiplied by the cut-off function $\chi_{[-y_0,0]}(y)$, such that
\begin{eqnarray}
\nonumber \frac{d {\bf c}_h}{d y} - \epsilon \Lambda_h({\bf c}_0)
{\bf c}_h & = & \epsilon \chi_{[0,y_0]}(y) {\bf G}_Q({\bf
c}_0)({\bf c}_h,{\bf c}) + \epsilon^{N+1} \chi_{[0,y_0]}(y) {\bf
G}_R({\bf c}_0 + {\bf c}_h,{\bf c}) \\ \label{semilinear-7} &
\phantom{t} & \epsilon \chi_{[-y_0,0]}(y) {\bf G}^*_Q({\bf
c}_0)({\bf c}_h,{\bf c}) + \epsilon^{N+1} \chi_{[-y_0,0]}(y) {\bf
G}^*_R({\bf c}_0 + {\bf c}_h,{\bf c}),
\end{eqnarray}
where $({\bf G}^*_{Q,R})_{1,2}(y) = (\bar{\bf G}_{Q,R})_{1,2}(-y)$
and $({\bf G}^*_{Q,R})_{3,4}(y) = -(\bar{\bf G}^*_{Q,R})_{1,2}(y)$
on $y \in [-y_0,0)$. We are looking for a global solution of the
system (\ref{semilinear-6})--(\ref{semilinear-7}) in the space of
bounded continuous functions $C_b^0(\mathbb{R})$. This global
solution on $y \in \mathbb{R}$ corresponds to a local solution of
the system (\ref{semilinear-4})--(\ref{semilinear-5}) on the
interval $y \in [0,y_0] \subset \mathbb{R}$.

\noindent The imaginary axis lies in the resolvent set of
$\Lambda_u$ and $\Lambda_s$ and
\begin{equation}
\label{stable-unstable} \| (\Lambda_{u,s} - i k I)^{-1}
\|_{X_{u,s} \mapsto X_{u,s}} \leq \frac{K_0}{1 + |k|}, \quad
\forall k \in \mathbb{R},
\end{equation}
for some $K_0 > 0$. Let ${\bf c}_s(y) = e^{y \Lambda_s} {\bf b} +
\tilde{\bf c}_s(y)$ and look for solution $\tilde{\bf c}_s(y)$ and
${\bf c}_u(y)$ of the system (\ref{semilinear-6}) projected to
$X_s$ and $X_u$ with operators $P_s$ and $P_u$. By Lemmas
\ref{lemma-inversion}, \ref{lemma-nonlinear-vector-field}, and the
Implicit Function Theorem, there exists a unique map from
$C^0_b(\mathbb{R},X_h \oplus X_c)$ to $C^0_b(\mathbb{R},X_u \oplus
X_s)$ parameterized by ${\bf b}$ such that
\begin{eqnarray}
\nonumber \sup_{y \in [0,y_0]} \left[ \| {\bf c}_u(y) \|_{X_u} +
\| {\bf c}_s(y) \|_{X_s} \right] & \leq & \| {\bf b} \|_{X_s} +
\epsilon M_1 \sup_{y \in [0,y_0]} \left[ \left( 1 + \| {\bf
c}_h(y) \|_{X_h} + \| {\bf c}_c(y) \|_{X_c} \right) \| {\bf c}_c(y) \|_{X_c} \right] \\
\label{mapping-stable-unstable} & \phantom{t} & + \epsilon^{N+1}
M_2 \sup_{y \in [0,y_0]} \left[1 + \| {\bf c}_h(y) \|_{X_h} + \|
{\bf c}_c(y) \|_{X_c} \right],
\end{eqnarray}
for some $M_1,M_2 > 0$.

\noindent Let ${\bf c}_h(y) = \alpha_1 {\bf s}_1(y) + \alpha_2
{\bf s}_2(y) + \tilde{\bf c}_h(y)$ and look for solution
$\tilde{\bf c}_h(y)$ of the system (\ref{semilinear-7}). By the
Hamiltonian structure of the system (\ref{semilinear-4}), the
vector field satisfies the constraints  $({\bf G}_{Q,R})_{3,4}(y)
= -(\bar{\bf G}_{Q,R})_{1,2}(y)$ on $y \in [0,y_0]$. By the
construction of the modified vector field, it satisfies the
constraint (\ref{constraints-kernel}). By Lemmas
\ref{lemma-weak-stable-unstable-manifolds},
\ref{lemma-nonlinear-vector-field}, the bound
(\ref{mapping-stable-unstable}), and the Implicit Function
Theorem, there exists a unique map from $C^0_b(\mathbb{R},X_c)$ to
$C^0_b(\mathbb{R},X_h)$ parameterized by $(\alpha_1,\alpha_2)$ and
${\bf b}$ such that
\begin{equation}
\label{mapping-stable-unstable-weak} \sup_{y \in [0,y_0]} \| {\bf
c}_h(y) \|_{X_h} \leq |\alpha_1| + |\alpha_2| + M_3 \sup_{y \in
[0,y_0]} \| {\bf c}_c(y) \|^2_{X_c} + \epsilon^N M_4 \sup_{y \in
[0,y_0]} \left[ 1 + \| {\bf c}_c(y) \|_{X_c} \right],
\end{equation}
for some $M_3,M_4 > 0$.

\noindent Since the spectrum of $\Lambda_c$ consists of pairs of
semi-simple purely imaginary eigenvalues, the operator $\Lambda_c$
generates a strongly continuous group $e^{y \Lambda_c}$ for any $y
\in \mathbb{R}$ on $X_c$ such that
\begin{equation}
\label{center} \sup_{y \in \mathbb{R}} \| e^{y \Lambda_c} \|_{X_c
\mapsto X_c} \leq K,
\end{equation}
for some $K > 0$. By variation of constant formula, the solution
of the system (\ref{semilinear-6}) projected to $X_c$ can be
rewritten in the integral form
\begin{equation}
\label{mapping-center} {\bf c}_c(y) = e^{y \Lambda_c} {\bf a} +
\epsilon \int_0^y  e^{(y-y') \Lambda_c} P_c \left[ {\bf F}_Q({\bf
c}_0(\epsilon y') + {\bf c}_h(y'),{\bf c}(y')) + \epsilon^N {\bf
F}_R({\bf c}_0(\epsilon y') + {\bf c}_h(y'),{\bf c}(y')) \right]
dy',
\end{equation}
where ${\bf a} = {\bf c}_c(0)$. By using the bound
(\ref{bound-homoclinic-orbit}) for ${\bf c}_0(\epsilon y)$ and the
bounds (\ref{mapping-stable-unstable}) and
(\ref{mapping-stable-unstable-weak}) on the components ${\bf
c}_{u,s}$ and ${\bf c}_h$, we derive from the integral equation
(\ref{mapping-center}) that
\begin{eqnarray*}
\sup_{y \in [0,y_0]} \| {\bf c}_c(y) \|_{X_c} \leq K \left( \|
{\bf a} \|_{X_c} + |{\bf b}\|_{X_s} + |\alpha_1| + |\alpha_2| +
\epsilon M_5 \int_0^{y_0} \| {\bf c}_0(y) \|_{X_h} \| {\bf c}_c(y)
\|_{X_c} dy \right. \\
\left. + \epsilon y_0 M_6 \sup_{y \in [0,y_0]} \| {\bf c}_c(y)
\|^2_{X_c} + \epsilon^{N+1} M_7 \int_0^{y_0} \| {\bf c}_0(y)
\|_{X_h} dy + \epsilon^{N+1} y_0 M_8 \sup_{y \in [0,y_0]} \| {\bf
c}_c(y) \|_{X_c} \right)
\end{eqnarray*}
for some $M_5,M_6,M_7,M_8 > 0$. By the Gronwall's inequality, we
have thus obtained that
\begin{eqnarray}
\nonumber \sup_{y \in [0,y_0]} \| {\bf c}_c(y) \|_{X_c} & \leq & K
e^{\epsilon K M_5 \int_0^{y_0} \| {\bf c}_0(y) \|_{X_h} dy} \;
\left( \| {\bf
a} \|_{X_c} + |{\bf b}\|_{X_s} + |\alpha_1| + |\alpha_2| \right. \\
\label{estimation-center} & \phantom{t} &  \left. + \epsilon^N M_9
+ \epsilon y_0 M_6 \sup_{y \in [0,y_0]} \| {\bf c}_c(y) \|^2_{X_c}
+ \epsilon^{N+1} y_0 M_8 \sup_{y \in [0,y_0]} \| {\bf c}_c(y)
\|_{X_c} \right),
\end{eqnarray}
for some $M_9 > 0$. Here we can use the decay bound
(\ref{bound-homoclinic-orbit}) which implies that $\epsilon
\int_0^{y_0} \| {\bf c}_0(y) \|_{X_h} dy  \leq C$ for some
$\epsilon$-independent $C > 0$. By using the same bound for the
exponent and letting $y_0 = L/\epsilon^M$, we can see that we can
choose $M \leq N+1$, where the value $M = N+1$ gives the balance
of all terms in the upper bound (\ref{estimation-center}). If
arbitrary vectors ${\bf a}$, ${\bf b}$ and $(\alpha_1,\alpha_2)$
satisfies the bound (\ref{bound-initial}), then we have
constructed a local solution ${\bf c}_c(y)$ which satisfies the
bound
\begin{equation}
\label{bound-final} \sup_{y \in [0,L/\epsilon^{N+1}]} \| {\bf
c}_c(y) \|_{X_c} \leq \tilde{C}_c \epsilon^N
\end{equation}
for some $\tilde{C}_c > 0$. By using the bounds
(\ref{bound-initial}), (\ref{mapping-stable-unstable}),
(\ref{mapping-stable-unstable-weak}), and (\ref{bound-final}), we
have proved the bound (\ref{bound-center}) for some
$\epsilon$-independent constants $C_h,C > 0$.
\end{proof}

\begin{remark}
{\rm One can prove Theorem \ref{theorem-reduction} by using the
contraction mapping principle and the integral formulation for the
local center manifold of the system
(\ref{semilinear-4})--(\ref{semilinear-5}). This approach was
undertaken in Section 4 of \cite{GS01} (see their Theorem 4). We
have avoided this unnecessary complication with the explicit
decomposition (\ref{decomposition-final}) and analysis of the
system (\ref{semilinear-4})--(\ref{semilinear-5}) decomposed into
subsystems. Similar direct methods of analysis have been applied
to problems without Hamiltonian structures such as the
quasi-linear wave equation in \cite{GS05}, where an iteration
scheme was employed to prove the bound on small local solutions
along the local center-stable manifold.}
\end{remark}

\begin{proof1}{\em of Theorem \ref{theorem-main}:} By Theorem
\ref{theorem-reduction}, we have constructed an
infinite-dimensional continuous family of local bounded solutions
of the system (\ref{semilinear-4})--(\ref{semilinear-5}) on $y \in
[0,L/\epsilon^{N+1}]$ for some $\epsilon$-independent constant $L
> 0$. The solutions are close to the reversible
homoclinic orbit of the extended coupled-mode system
(\ref{cme-normal-form}) in the sense of the bound
(\ref{bound-main}). It remains to extend the local solution to the
symmetric interval $y \in [-L/\epsilon^{N+1},L/\epsilon^{N+1}]$ as
the local reversible solution with the reversibility constaints
(\ref{domain-reversibility-modified}). To do so, we shall consider
the intersections of the local invariant manifold of the system
(\ref{semilinear-4})--(\ref{semilinear-5}) with the symmetric
section $\Sigma_r'$ defined by
(\ref{constraint-reversibility-modified}).

\noindent Since the initial data ${\bf c}_c(0) = {\bf a}$ in the
local center--stable manifold of Theorem \ref{theorem-reduction}
are arbitrary, the components of ${\bf a}$ can be chosen to lie in
the symmetric section $\Sigma_r'$, such that
\begin{equation}
{\rm Im} ({\bf a})_m^+ = 0, \; \forall m \in
\mathbb{Z}_-'\backslash \{n\}, \qquad {\rm Im} ({\bf a})_m^- = 0,
\; \forall m \in \mathbb{Z}_-'\backslash \{-n\}.
\end{equation}
This construction still leaves infinitely many arbitrary
parameters for
\begin{equation}
{\rm Re} ({\bf a})_m^+, \; \forall m \in \mathbb{Z}_-'\backslash
\{n\}, \qquad {\rm Re} ({\bf a})_m^-, \; \forall m \in
\mathbb{Z}_-'\backslash \{-n\}
\end{equation}
to be chosen in the bound (\ref{bound-initial}). The initial data
${\bf c}_h(0)$ and ${\bf c}_{s,u}(0)$ are not arbitrary since we
have used the Implicit Function Theorem for the mappings
(\ref{mapping-stable-unstable}) and
(\ref{mapping-stable-unstable-weak}). Therefore, we have to show
that the components of ${\bf b}$ and $(\alpha_1,\alpha_2)$ can be
chosen uniquely so that the local center-stable manifold
intersects at $y = 0$ with the symmetric section $\Sigma_r'$.

\noindent We note that there are as many arbitrary parameters
${\bf s}$ and $(\alpha_1,\alpha_2)$ in the local center--stable
manifold as there are remaining constraints in the set
$\Sigma_r'$. First, let us consider constraints in the set
$\Sigma_r'$ for $m \in \mathbb{Z}'_+$, namely
\begin{equation}
{\rm Re} c_m^+(0) = {\rm Re} c_m^-(0), \quad {\rm Im} c_m^+(0) =
-{\rm Im} c_m^-(0), \quad m \in \mathbb{Z}'_+.
\end{equation}
Let ${\bf c}_s = e^{y \Lambda_s} {\bf b} + \tilde{\bf c}_s(y)$ and
rewrite the constraints in the form
\begin{equation}
{\rm Re} b_m + {\rm Re}(\tilde{\bf c}_s)_m(0) = {\rm Re}({\bf
c}_u)_m(0), \quad {\rm Im} b_m + {\rm Im}(\tilde{\bf c}_s)_m(0) =
-{\rm Im}({\bf c}_u)_m(0),
\end{equation}
where all terms are of order ${\rm O}(\epsilon^N)$ and the vectors
$\tilde{\bf c}_s(0)$ and ${\bf c}_u(0)$ depend on ${\bf b}$ in
higher orders in $\epsilon$. By the Implicit Function Theorem,
there exists a unique solution of the constraints for ${\bf b}$
such that $\| {\bf b} \|_{X_s}$ satisfies the bound
(\ref{bound-initial}).

\noindent Finally, let us consider constraints in the set
$\Sigma_r'$ for components of ${\bf c}_h$, namely
\begin{equation}
\label{last-set-constraints} {\rm Im} c_n^+(0) = 0, \quad {\rm Im}
c_{-n}^-(0) = 0.
\end{equation}
Let ${\bf c}_h = \alpha_1 {\bf s}_1 + \alpha_2 {\bf s}_2 +
\tilde{\bf c}_s(y)$ and note that ${\bf s}_1$ and ${\bf s}_2$
violate the constraints (\ref{last-set-constraints}). Let ${\bf
a}$ and ${\bf b}$ be chosen so that ${\bf c}(0)$ belongs to the
set $\Sigma_r'$. By a construction of the vector field in the
system (\ref{semilinear-7}), if ${\bf c}(y)$ lies in the domain
$D_r'$ of reversible solution and $\alpha_1 = \alpha_2 = 0$, then
${\bf G}_{Q,R}^* = {\bf G}_{Q,R}$ on $y \in [-y_0,y_0]$ and the
global solution $\tilde{\bf c}_s(y)$ constructed in Theorem
\ref{theorem-reduction} intersects the set $\Sigma_r'$ at $y = 0$.
Therefore, the choice $\alpha_1 = \alpha_2 = 0$ satisfies the
constraints (\ref{last-set-constraints}) identically.

\noindent We have thus constructed a family of reversible
solutions in the symmetric interval $y \in
[-L/\epsilon^{N+1},L/\epsilon^{N+1}]$ while preserving the bound
(\ref{bound-main}). Tracing the coordinate transformations used in
our analysis back to the original variable $\psi(x,y)$, we have
thus completed the proof of Theorem \ref{theorem-main}.
\end{proof1}

\section{Discussion}

We have proved that a moving gap soliton of the Gross--Pitaevskii
equation (\ref{GP}) with the periodic potential $V(x)$ is
surrounded by the oscillatory tails which are bounded on {\em
finite} intervals of the spatial scale. Because the center
manifold is infinite-dimensional with the sign-indefinite
Hamiltonian function, we are not able to exclude the polynomial
growth of the oscillatory tails in the far-field regions. This
construction of traveling solutions on a finite spatial scale is
related with the {\em finite-time} applicability of the
coupled-mode equations (\ref{cme}) for the Cauchy problem
associated with the Gross--Pitaevskii equation (\ref{GP})
\cite{SU}.

It would have been a drastic improvement to the constructed theory
if we could extend the analysis of oscillatory tails to the {\em
infinite} spatial scale by proving existence of {\em global}
solutions with a single bump and {\em bounded} oscillatory tails.
In many problems with finite-dimensional center manifolds
associated with semi-simple purely imaginary eigenvalues, such
constructions of global center-stable manifolds are possible and
the proof of persistence of bounded solutions with oscillatory
tails can be developed \cite{IL}.

We will show that the basic evolution models for moving gap
solitons in periodic potentials exhibit infinite-dimensional
center manifolds in the spatial dynamics formulation. In
particular, we can think of three possible generalizations of the
Gross--Pitaevskii equation (\ref{GP}), given by the complex-valued
Klein--Gordon equation
\begin{equation}
\label{KG} E_{tt} - E_{xx} = V(x) E + \sigma |E|^2 E,
\end{equation}
the regularized Gross--Pitaevskii equation
\begin{equation}
\label{regGP} i E_t = - E_{xx} + i E_{xxt} + V(x) E + \sigma |E|^2
E,
\end{equation}
and the discrete Gross--Pitaevskii equation
\begin{equation}
\label{disGP} i \dot{E}_n = -E_{n+1}-E_{n-1} + V_n E_n + \sigma
|E_n|^2 E_n.
\end{equation}
When $V(x) \equiv 0$ or $V_n \equiv 0$, the spectrum of the linear
part of the Klein--Gordon equation (\ref{KG}) is unbounded from
both above and below, while that of the regularized and discrete
Gross--Pitaevskii equations (\ref{regGP}) and (\ref{disGP}) is
bounded from both above and below. We look at the traveling
solutions of these equations in the form
\begin{equation}
\label{traveling-wave-final} E(x,t) = \sum_{m \in \mathbb{Z}'}
\psi_m(y) e^{\frac{i m x}{2} - i \omega t}, \quad y = x - ct,
\end{equation}
for the linear limit $\sigma = 0$ with no potential $V(x) \equiv
0$. (In the case of the lattice equation (\ref{disGP}), we use the
traveling ansatz (\ref{traveling-wave-final}) at $x = n$ for $n
\in \mathbb{N}$.) As a result, we obtain uncoupled linear ODEs or
differential advance-delay equations for amplitudes $\{ \psi_m(y)
\}_{m \in \mathbb{Z}'}$ which are solved with the substitution
$\psi_m(y) = e^{\kappa y}$, $\forall m \in \mathbb{Z}'$. All roots
$\kappa$ are found from the following characteristic equations
\begin{eqnarray}
\label{ch-eq-KG} (1 - c^2) \kappa^2 + i (m - 2
c \omega) \kappa + \omega - \frac{m^2}{4} & = & 0, \\
\label{ch-eq-regGP} - i c \kappa^3 + (\omega
- 1 + m c) \kappa^2 + i \left( c - m + m \omega + \frac{m^2 c}{4} \right)
\kappa - \omega + \frac{m^2}{4} (1-\omega) & = & 0, \\
\label{ch-eq-disGP} \omega - i c \kappa + 2 \cosh \left(\kappa +
\frac{m}{2}\right) & = & 0,
\end{eqnarray}
which correspond to the three relevant models
(\ref{KG})--(\ref{disGP}). It is easy to see that the
characteristic equations (\ref{ch-eq-regGP}) and
(\ref{ch-eq-disGP}) have at least one purely imaginary root $p$
for any $m \in \mathbb{Z}'$, while the characteristic equation
(\ref{ch-eq-KG}) has two purely imaginary roots $p$ for
sufficiently large values of $|m|$ on $m \in \mathbb{Z}'$.
Therefore, the dimension of the center manifold associated with
the linearized system at the zero solution for $V(x) \equiv 0$ is
infinite in all three models (\ref{KG})--(\ref{disGP}).

{\bf Acknowledgement.} The work of D. Pelinovsky is supported by
the Humboldt Research Foundation. The work of G. Schneider is
partially supported  by the Graduiertenkolleg 1294 ``Analysis,
simulation and design of nano-technological processes'' granted by
the Deutsche Forschungsgemeinschaft (DFG) and the Land
Baden-W\"{u}rttemberg.

\end{document}